\documentclass[letterpaper, 11pt]{article}
\pdfoutput=1
\usepackage[letterpaper, margin=1in]{geometry}

\usepackage{bm}
\usepackage{amsmath}
\usepackage{subfigure}
\usepackage{amssymb}
\usepackage{mathtools}	
\usepackage[colorlinks=true, pdfstartview=FitV, linkcolor=black, citecolor= black, urlcolor= black]{hyperref}
\usepackage{overcite}
\usepackage{footnpag}			      	% make footnote symbols restart on each page
\usepackage{lipsum}
\usepackage{caption}
\usepackage{algorithm}% http://ctan.org/pkg/algorithm
\usepackage{algpseudocode}% http://ctan.org/pkg/algorithmicx

\usepackage{graphicx}
\usepackage{bm,amsmath,xfrac} 
\usepackage[version=4]{mhchem}
\usepackage[mathscr]{eucal}
\usepackage{siunitx}
\usepackage{longtable,tabularx}
\usepackage{nicematrix}

\usepackage{smartdiagram}
\usesmartdiagramlibrary{additions}

\newlength\figureheight 
\newlength\figurewidth
\newlength\insetwidth

\usepackage{pgfplots,tikz,tikz-3dplot}
\usepackage{tikzscale}
\usetikzlibrary{shapes,backgrounds,plotmarks, arrows}
\pgfplotsset{compat=newest} 
\pgfplotsset{plot coordinates/math parser=false}
\usepackage[10pt]{moresize}
\pgfplotsset{every tick label/.append style={font=\tiny}}

\newcommand{\equ}[1]{Eq.~(\ref{#1})}
\newcommand{\equt}[2]{Eqs.~(\ref{#1})~and~(\ref{#2})}

\newcommand{\bequ}[1]{Equation~(\ref{#1})}
\newcommand{\bequt}[2]{Equations~(\ref{#1})~and~(\ref{#2})}
\newcommand{\bequs}[2]{Equations~(\ref{#1})-(\ref{#2})}
\newcommand{\bequg}[1]{Equations~(\ref{#1})}

\newcommand{\bfig}[1]{Figure~(\ref{#1})}
\newcommand{\bfigt}[2]{Figures~(\ref{#1})~and~(\ref{#2})}
\newcommand{\bfigs}[2]{Figures~(\ref{#1})-(\ref{#2})}

\newcommand{\tbl}[1]{Table~\ref{#1}}

\newcommand{\exv}[1]{\operatorname{\mathbb{E}}\lbrace#1\rbrace}

\newcommand{\td}[0]{\text{d}}

\renewcommand{\l}[0]{(\ell)}

\begin{document}

\title{\textbf{Probabilistic Trajectory Design Via Approximate Gaussian Mixture Steering}}

\date{}

\author{William N. Fife\thanks{Graduate Student, Department of Aerospace Engineering, Texas A\&M University, College Station, TX.}, 
Pradipto Ghosh\thanks{Sr. Professional Staff II , SES/SAC , Johns Hopkins University Applied Physics Laboratory, 11100 Johns Hopkins Road, Laurel, MD 20723-6099.}, 
and Kyle J. DeMars\thanks{Associate Professor, Department of Aerospace Engineering, Texas A\&M University, College Station, TX.}
\ 
%\ and XXXXXX\thanks{YYYYY, YYYYY, YYYYY, YYYYY YYYYY.}
}

\maketitle{}

\begin{abstract}
A method is presented to solve a stochastic, nonlinear optimal control problem representative of spacecraft trajectory design under uncertainty. The problem is reformulated as a chance constrained nonlinear program, or what is known as a distribution steering problem. Typical distribution steering problems rely on the underlying uncertainties to be Gaussian distributions. This work expands on previous developments by embedding Gaussian mixture distributions into the formulation to better handle the uncertainty propagation and chance constraints involved. The method is applied to a finite-thrust Earth-to-Mars transfer problem. Evaluation via Monte Carlo analysis shows a greater satisfaction of constraints under non-Gaussian distributions of the state and a statistically lower cost. 
\end{abstract}

\section*{Introduction}
The task of designing a spacecraft's trajectory, whether to a celestial body or another spacecraft, is typically solved using optimal control theory\cite{bryson2018applied}. In all but trivial cases, the optimal control problem is converted into a parameter optimization problem\cite{hull1997conversion} that is computationally solved to yield a set of nominal controls (inputs) that steer the spacecraft state to a desired condition by minimizing a cost function. Typically, these methods rely on perfect modeling of the dynamical evolution of the state and/or perfect knowledge of the state. To account for uncertainty, it is usual to augment the nominal sequence of controls with a feedback strategy based on results from Monte Carlo (MC) analyses that simulate truth models that are different from the deterministic models used in the original optimization. This process is not only computationally expensive and lengthy, but also must be finely tuned for the specific problem data. Alternatively, a design process, such as the one employed in this work, that embeds the model and state uncertainties directly into the optimization of the controls ameliorates these shortcomings.

Model and state knowledge imperfections are traditionally dealt with by considering the system to be stochastic. It has been shown\cite{hotz1987covariance, bakolas2016optimal}, for controllable linear systems, that an initial state covariance can be steered by a feedback law to some prescribed covariance while minimizing control effort. This method is referred to as \textit{covariance steering}. For spacecraft trajectory design, there are multiple difficulties associated with implementing traditional covariance steering methods. Firstly, spacecraft dynamics are generally nonlinear and thus pose modeling difficulties \textit{a priori} of stochasticity. Additionally, minimization of control effort in a deterministic sense may not guarantee that the control at any given instant does not exceed a magnitude physically unrealizable by the spacecraft's propulsion system. Therefore, constraints on the control must be included in the stochastic optimization problem. Since the control is formed by a feedback law, the control is a stochastic quantity, and violation of constraints on the control are uncertain. This is dealt with by considering \textit{chance constraints} that bound probabilities of control magnitude constraints.

For Gaussian systems, the cumulative distribution function and other concentration inequalities\cite{wainwright2019high} allow one to convert chance constraints into equivalent deterministic ones. Unfortunately, the nonlinear spacecraft dynamics propagated over typical interplanetary mission durations can create significantly non-Gaussian distributions. Approximating nonlinear and/or non-Gaussian chance constraints is an active research topic in many fields (e.g., engineering, finance, etc.). Previous works have used linearization\cite{nemirovski2012safe}, sampling\cite{pagnoncelli2009sample}, or parametric\cite{geletu2015tractable} methods to approximate these chance constraints. This work seeks to extend the state-of-the-art in covariance steering by developing an approximate Gaussian mixture (GM) nonlinear program that is well-suited to handle the nonlinearity of the dynamics and non-Gaussianity of the chance constraints. GM-based steering (appropriately titled \textit{density steering}) has been recently investigated as a way to rigorously lift the original problem out of linear-Gaussian assumptions\cite{balci2023density}. For example, the work by Boone and McMahon\cite{boone2022spacecraft} uses a GM relaxation of chance constraints by considering weighted chance constraints for each component in the GM. Although rigorous, evaluating constraints on each GM component can become unwieldy for mixtures with a large number of components that are sometimes necessary to faithfully capture the evolution of uncertainty. Further, the weights for each component chance constraint need to be optimized to reduce unnecessary conservatism. 

The procedure developed in this work is for continuous/discrete dynamical systems whereby the control is applied at discrete epochs (nodes). A given Gaussian distribution (or a set of particles), representing the initial state uncertainty, is accurately approximated by a GM and propagated through the continuous, nonlinear dynamics via linearization. At each node, the GM is collapsed to a single Gaussian to evaluate a Gaussian chance constraint. The process is then repeated for each discrete node in the optimization. This approach more accurately models the true evolution of the state uncertainty compared to those that propagate a single Gaussian over the optimization. The application of this method reveals a lower statistical cost and tighter satisfaction of the terminal constraints than those reported elsewhere in the literature for the same test problem\cite{marmo2022chance}.

This paper is organized as follows. First, the relevant stochastic optimal control problem is described along with traditional assumptions made. A Gaussian mixture uncertainty framework is then derived for approximating and propagating the true distribution through the nonlinear dynamics. Then the presented optimal control problem, along with the converted chance constraints, are summarized into a nonlinear program. Lastly, the presented algorithm is applied to an Earth-to-Mars finite thrust scenario and evaluated via Monte Carlo analysis.

\section*{Problem Description}
The stochastic state $\bm{x} \in \mathbb{R}^n$ considered in this work is governed by hybrid dynamics of the form
\begin{subequations}
    \label{eqn::dyn}
    \begin{align}
        \label{eqn::dyna}
        \dot{\bm{x}}(t) &= \bm{f}(\bm{x}(t)) + \bm{F}_w(t)\bm{w}(t) \qquad t \in [t_{k-1}, t_k] \\
        \label{eqn::dynb}
        \bm{x}_{k_+} &= \bm{x}_{k_-} + \bm{F}_u\bm{u}_k \, ,
    \end{align}
\end{subequations}
 where $\bm{w}(t) \in \mathbb{R}^p$ is a zero-mean, process noise term with covariance $\bm{P}_{ww}(t) = \bm{Q}(t)\delta(t - \tau)$, with $\bm{Q}(t)$ being a time-varying power spectral density. The subscripts $k_-$ and $k_+$ denote a prior and posterior quantity to the node at time $t_k$, respectively. The control $\bm{u}_k \in \mathbb{R}^m$ is applied at discrete nodes via the control mapping matrix $\bm{F}_u$. The process noise is assumed to be uncorrelated with the state at the initial time. At the initial time $t = t_0$, the set of all possible states is characterized by some general probability density function (pdf) $p(\bm{x}_0)$. The objective is to compute a set of feedback controls $\bm{u}_{0:N-1} = \lbrace \bm{u}_0(\bm{x}_0), \bm{u}_1(\bm{x}_1), ..., \bm{u}_{N-1}(\bm{x}_{N-1}) \rbrace$ that steer the state from $p(\bm{x}_0)$ to some prescribed pdf $p(\bm{x}_N)$ in a given finite interval $t \in [t_0, t_N]$ by minimizing a scalar cost function $J$ and adhering path constraints on the control and state.

In this work, the boundary constraints on the state are the aforementioned pdf prescriptions, and the path constraints are given by the dynamics in \bequg{eqn::dyn}. The control has no boundary constraints, but has chance path constraints, applied node-wise, given by the probability evaluation
\begin{align}
    \label{eqn::prob-u}
    P(\| \bm{u}_k \| \leq \rho_u) \geq 1 - \beta \qquad \forall k \in [0, N-1] \, ,
\end{align}
where $\beta$ is a small, positive number that represents the prescribed probability of the constraint being violated and $\rho_u$ is a prescribed maximum control norm. The motivation behind the bounded control constraint being probabilistic is to rigorously account for uncertainties within the control itself. These uncertainties stem from multiple sources including, but not limited to, random noise, biases, and state feedback uncertainties. The cost function considered in this work is given by
\begin{align}
    \label{eqn::cost}
    J = \sum_{k = 0}^{N-1} \| \exv{\bm{u}_k} \| \, ,
\end{align}
which is typically used to mimic minimum control objectives in stochastic systems \cite{ross2004find}. Therefore, the nonlinear, optimal uncertainty steering problem is summarized as
\begin{subequations}
    \label{eqn::USprob}
    \begin{align}
        \label{eqn::USprob-a}
        \min_{\bm{u}_{0:N-1}} \sum_{k = 0}^{N-1} \| \exv{\bm{u}_k} \| \, , \quad \qquad& \\
        \text{subject to path constraints}& \nonumber \\
        \label{eqn::USprob-b}
        \dot{\bm{x}}(t) &= \bm{f}(\bm{x}(t)) + \bm{F}_w(t)\bm{w}(t) \qquad &t \in [t_{k-1}, t_k] \\
        \label{eqn::USprob-c}
        \bm{x}_{k_+} &= \bm{x}_{k_-} + \bm{F}_u\bm{u}_k \qquad &\forall k \in [0, N-1] \\
        \label{eqn::USprob-d}
        P(\| \bm{u}_k \| \leq \rho_u) &\geq 1 - \beta \qquad &\forall k \in [0, N-1] \\
        \text{and terminal constraints}& \nonumber \\
        \label{eqn::USprob-e}
        \bm{x}_0 &\sim p(\bm{x}_0) \\
        \label{eqn::USprob-f}
        \bm{x}_N &\sim g(\bm{x}_N) \, .
    \end{align}
\end{subequations}
\bequg{eqn::USprob} represent a generic architecture for the minimum-control uncertainty steering problem subject to continuous-discrete dynamics. Unfortunately, there is no general solution to the problem, but only a framework for which satisfying (exactly or approximately) certain conditions yields a potential solution. The three primary conditions are: i) for a given control parameterization, \bequs{eqn::USprob-b}{eqn::USprob-c} can be decomposed into equations governing the moments of the pdf of $\bm{x}$, ii) that those moments are sufficient statistics for both $p(\bm{x}_0)$ and $p(\bm{x}_N)$, and iii) \bequs{eqn::USprob-e}{eqn::USprob-f} can be decomposed into deterministic constraints on the moments that are feasible. If those conditions are exactly met, then \bequg{eqn::USprob} can be converted into a nonlinear programming problem that, if solved, yields a control that guarantees the initial pdf be optimally\footnote{Optimum is not guaranteed to be global unless the problem is convex.} steered to the final pdf.

Typically, it is assumed that the initial pdf is a multivariate Gaussian and that the control is a linear feedback policy (either Markovian or with history) \cite{skaf2010design, bertsekas1978stochastic}. The complication is that the presence of nonlinear state dynamics means there is no guarantee the state pdf remains Gaussian. Approximate formulations are created by forming a Taylor series approximation to the dynamics or using some form of sampling integration \cite{sarkka2023bayesian}. Over long integration periods, these approximations can provide unsatisfactory results. Moreover, a Gaussian assumption of the state pdf restricts the study of more exotic problems with multi-modal or heavy-tailed distributions. Therefore, it is desirable to use less restrictive uncertainty representations within the problem framework. 

\section*{Methodology}
This section outlines the decomposition of the nonlinear, optimal uncertainty steering problem summarized in \bequg{eqn::USprob}. The Gaussian mixture framework for uncertainty propagation is developed along with a model for incorporating maneuver execution error. Then, the cost and constraint functions are formulated in terms of the optimization variables. Lastly, the developments are summarized into a nonlinear program. 

\subsection*{Gaussian Mixture Uncertainty Framework}
\begin{figure}[htbp]
    \centering
    \setlength\figurewidth{8cm}
    \setlength\figureheight{4.5cm}
    \input{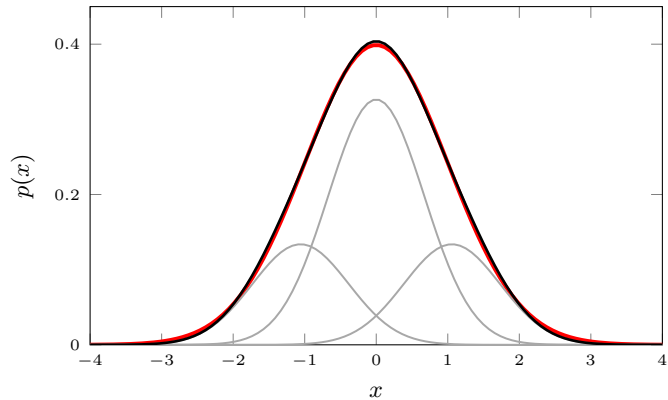}
    \caption{GM approximation of standard normal using $\bm{L = 3}$ splitting library. The standard Gaussian is given by the red line, the individual GM components are shown in gray, and the overall GM is shown in black.}
    \label{fig::L3GM}
\end{figure}

The initialization of the state pdf can take various forms depending on the \textit{a priori} information on the state uncertainty. No initialization procedures are required if already starting with a GM. Many orbit determination algorithms in use today return a single estimate (mean) and covariance representing a Gaussian uncertainty, so that is the starting point in this work. Alternatively, an expectation maximization algorithm can be used to initialize a GM from a set of particles\cite{mclachlan2019finite}. Let the state pdf at initial time $t_0$ be given by a Gaussian distribution with known mean $\bm{m}_{x,0}$ and covariance $\bm{P}_{xx,0}$. This single multivariate Gaussian distribution is then split into a multivariate GM pdf that closely approximates the original Gaussian distribution. There are several algorithms available to complete this step \cite{li2009novel, hanebeck2003progressive}. The method used in this work was developed by DeMars \cite{demars2013entropy}. A brief outline of the algorithm is as follows. First, consider the univariate case of finding the weights $\Tilde{w}^{\l}$, means $\Tilde{m}^{\l}$, and variance $\Tilde{\sigma}^2$ of a GM pdf that approximates a standard Gaussian; that is,
\begin{align*}
    \sum_{\ell = 1}^{L} \Tilde{w}^{\l} p_g(x; \Tilde{m}^{\l}, \Tilde{\sigma}^2) \approx p_g(x; 0, 1) \, .
\end{align*}
It is noted that the GM approximation is constrained to have components with identical variance. Let the GM pdf be denoted by $\Tilde{p}(x)$ and the Gaussian by $p(x)$. A minimization problem is posed with a cost given by
\begin{align}
    \label{eqn::splitcost}
    J = \int_{\mathcal{S}} ( p(x) - \Tilde{p}(x) )^2 \td x + \lambda \Tilde{\sigma}^2 \quad \text{subject to} \quad \sum_{\ell=1}^{L} \Tilde{w}^{\l} = 1 \, ,
\end{align}
where the first term on the right-hand side of \bequ{eqn::splitcost} is the $L_2$ distance between the pdfs over common support $\mathcal{S}$, and $\lambda$ is a weighting term to scale the cost of larger $\Tilde{\sigma}^2$. For a given number of components $L$, this cost can be minimized and the result is a library of optimal weights, means, variances. This work uses libraries computed for $L = \lbrace3, 4, 5\rbrace$. An example library for $L = 3$ and $\lambda = 0.001$ is shown in \tbl{tbl::split3} with an example given in \bfig{fig::L3GM}.
\begin{table}[htbp]
    \centering
    \caption{Splitting library for $\bm{L = 3}$ and $\bm{\lambda = 0.001}$.}
    \label{tbl::split3}
    \begin{NiceTabular}{cccc} \hline \hline
        $\ell$ & $\tilde{w}^{\l}$ & $\tilde{m}^{\l}$ & $\tilde{\sigma}$ \\
        \hline $1$ & $0.2252246249$ & $-1.0575154615$ & $0.6715662887$ \\
        $2$ & $0.5495507502$ & $0$ & $0.6715662887$ \\
        $3$ & $0.2252246249$ & $1.0575154615$ & $0.6715662887$ \\
        \hline \hline
    \end{NiceTabular}
\end{table}
The univariate splitting technique is applied to the multivariate Gaussian by computing the square root factor $\bm{S}_{xx}$ such that $\bm{S}_{xx}\bm{S}_{xx}^T = \bm{P}_{xx}$, and splitting along the directions given by the columns of $\bm{S}_{xx}$. This is done recursively for each dimension desired to be split, resulting in a total of $L = N^d$ components, where $N$ is the number of components being split into for each dimension, and $d$ is the number of dimensions desired to be split. The end result is a state pdf at initial time $t_0$ expressed as
\begin{align}
    \label{eqn::initpdf}
    p(\bm{x}_0) = \sum_{\ell = 1}^{L} w_{x,0}^{\l} p_g(\bm{x}_0 ; \bm{m}_{x,0}^{\l}, \bm{P}_{xx,0}^{\l}) \, ,
\end{align}
where the weights $w_x^{\l}$ sum to unity. It is important to note that the splitting library is pre-computed so that the actual splitting of the Gaussian is a fast procedure. 

\subsubsection*{Gaussian Mixture Propagation}
For the sake of brevity, the entire derivation for the propagation of the GM pdf through the continuous time dynamics is not included in this work. However, it is important to show the approximations used in the developments---the primary one being that the nonlinear dynamics are approximated by a first-order Taylor series expansion about the mean as
\begin{align}
    \label{eqn::FOTSE}
    \bm{f}(\bm{x}(t)) \approx \bm{f}(\bm{m}_x(t)) + \bm{F}_x(\bm{m}_x(t))[\bm{x}(t) - \bm{m}_x(t)] \, ,
\end{align}
where $\bm{F}_x(\bm{m}_x(t))$ is the Jacobian of $\bm{f}(\bm{x}(t))$ evaluated at the mean. It is also assumed that the process noise is Gaussian at all times. With these assumptions, it has been shown \cite{sarkka2019applied} that the continuous time evolution of the GM pdf is as follows. For each $\ell$-th component of the GM, the mean is propagated from $t_{k-1}$ to $t_{k_-}$ via
\begin{align}
    \label{eqn::meanprop}
    \dot{\bm{m}}_x^{(\ell)} (t) = \bm{f}(\bm{m}_x^{(\ell)}(t)) \, ,
\end{align}
with initial condition $\bm{m}_x^{(\ell)}(t_{k-1}) = \bm{m}_{x,k-1}^{(\ell)}$. Simultaneously, the $\ell$-th covariance is propagated by
\begin{align}
    \label{eqn::covprop}
    \dot{\bm{P}}_{xx}^{(\ell)} (t) &= \bm{F}_x(\bm{m}_x^{(\ell)}(t)) \bm{P}_{xx}^{(\ell)}(t) + \bm{P}_{xx}^{(\ell)}(t) \bm{F}_x^T(\bm{m}_x^{(\ell)}(t)) + \bm{F}_w(t)\bm{Q}(t)\bm{F}_w^T(t) \, ,
\end{align}
with initial condition $\bm{P}_{xx}^{(\ell)}(t_{k-1}) = \bm{P}_{xx,k-1}^{(\ell)}$. The component weights are held constant across the integration period, such that $w_{x,k_-}^{\l} = w_{x,k-1}^{\l}$. The number of components is also constant through propagation. The result is a GM pdf at time $t_{k_-}$ given by
\begin{align}
    \label{eqn::GMprop}
    p(\bm{x}_{k_-}) = \sum_{\ell = 1}^{L} w_{x,k_-}^{\l} p_g(\bm{x}_{k_-}; \bm{m}_{x,k_-}^{\l}, \bm{P}_{xx,k_-}^{\l}) \, .
\end{align}
Before the node update, at $t_{k_-}$, the GM pdf is collapsed down to a single mean and covariance such that the state pdf at $t_{k_-}$ is assumed to be $p_g(\bm{x}_{k_-}; \bm{m}_{x,k_-}, \bm{P}_{xx,k_-})$ where
\begin{align}
    \label{eqn::GMmean}
    \bm{m}_{x,k_-} = \sum_{\ell = 1}^{L}w_{x,k_-}^{\l} \bm{m}_{x,k_-}^{\l} 
\end{align}
and
\begin{align}
    \label{eqn::GMcov}
	\bm{P}_{xx,k_-}=\sum_{\ell=1}^{L} w_{x,k_-}^{(\ell)}\left(\boldsymbol{P}_{x x, k_-}^{(\ell)}+(\bm{m}_{x,k_-}^{\l}-\bm{m}_{x,k_-})(\bm{m}_{x,k_-}^{\l}-\bm{m}_{x,k_-})^T\right) \, .
\end{align}
The collapse of the GM pdf back into a single Gaussian represents an assumption that the mean and covariance encapsulate most of the probability mass of the true pdf of $\bm{x}$ at $t_{k-}$. The motivation behind the collapse is to approximate the control as a Gaussian random vector in order to use a previously-developed chance-constraint conversion. A future development of this work is envisioned where this assumption will be more rigorously addressed.

\subsubsection*{Control Parameterization and Node Update}

\begin{figure}[h!]
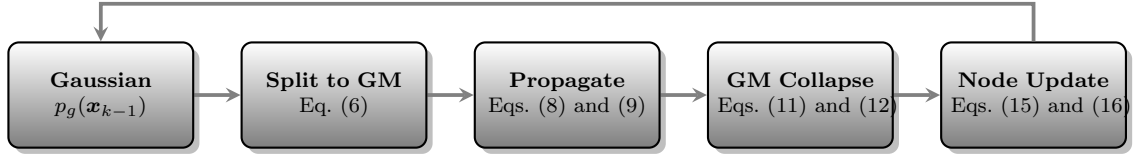

    \centering
    \smartdiagramset{uniform color list=gray for 5 items,
    module minimum width=2.2cm,
    module minimum height=1.4cm,
    module x sep=3.1,
    text width=2.2cm,
    font=\scriptsize,
    border color=black,
    arrow line width=0.05cm}
    \smartdiagram[flow diagram:horizontal]{\textbf{Gaussian} \\ $p_g(\bm{x}_{k-1})$,
        \textbf{Split to GM} \\ \equ{eqn::initpdf}, \textbf{Propagate} \\ \equt{eqn::meanprop}{eqn::covprop}, \textbf{GM Collapse} \\ \equt{eqn::GMmean}{eqn::GMcov}, \textbf{Node Update} \\ \equt{eqn::meanupdate}{eqn::covupdate}}
    \caption{Flow chart of propagation through dynamics and control update with associated equations. Recursion is established by $\bm{t_{k-1} \leftarrow t_{k_+}}$.}
    \label{fig::flowchart}
\end{figure}

To simplify the generic problem given by \bequg{eqn::USprob}, the control is parameterized by the affine Markov feedback law
\begin{align}
    \label{eqn::control}
    \bm{u}_k = \bm{v}_k + \bm{G}_k [\bm{x}_{k_-} - \bm{m}_{x,k_-}] + \delta \bm{u}_k \, ,
\end{align}
where $\bm{v}_k$ represents a nominal or open-loop control and $\bm{G}_k$ is a feedback gain maxtrix operating on the discrepancy between the actual/random state and the mean. The $\delta \bm{u}_k$ term acts as a random variable representing errors in the execution of the control. There are various ways to model this execution error, such as the Gates' model \cite{gates1963simplified}. In this work, $\delta \bm{u}_k$ is modeled as a zero-mean Gaussian random variable with prescribed covariance $\bm{P}_{\delta u \delta u, k}$. With \bequ{eqn::control}, the node update in \bequ{eqn::USprob-c} becomes
\begin{align}
    \label{eqn::nodeupdate}
    \bm{x}_{k_+} = \bm{x}_{k_-} + \bm{F}_u\left( \bm{v}_k + \bm{G}_k [\bm{x}_{k_-} - \bm{m}_{x,k_-}] + \delta \bm{u}_k \right) \, .
\end{align}
Assuming that $\exv{\bm{x}_{k_-}} = \bm{m}_{x,k_-}$, the mean of the node-updated state is given by
\begin{align}
    \label{eqn::meanupdate}
    \bm{m}_{x,k_+} = \bm{m}_{x,k_-} + \bm{F}_u \bm{v}_k \, .
\end{align}
Recognizing that $\bm{F}_u$ and $\bm{G}_{k}$ are deterministic quantities, computing $\exv{[\bm{x}_{k_+} - \bm{m}_{x,k_+}][\bm{x}_{k_+} - \bm{m}_{x,k_+}]^T}$ gives the covariance of the node-updated state as
\begin{align}
    \label{eqn::covupdate}
    \bm{P}_{xx,k_+} = [\bm{I}_n + \bm{F}_u \bm{G}_k]\bm{P}_{xx,k_-}[\bm{I}_n + \bm{F}_u \bm{G}_k]^T + \bm{P}_{\delta u \delta u, k} \, .
\end{align}
The final step of the node update is to establish the recursion $\bm{m}_{x,k-1} = \bm{m}_{x,k_+}$ and $\bm{P}_{xx,k-1} = \bm{P}_{xx,k_+}$. The entire procedure is then repeated for $N-1$ nodes in time. The flow chart in \bfig{fig::flowchart} highlights the primary steps with the associated equations. 
\vspace{1mm}

A note about the number of GM components is relevant for the developments herein: in this work, the number of components for GM splitting is kept the same for every iteration through the nodes. This is not a limitation of the current formulation. It is possible to implement any number of adaptive splitting methods during the continuous-time propagation \cite{demars2013entropy, sun2023hybrid}. For the analysis performed in this work, the authors found it sufficient not to implement an adaptive splitting method.

\subsection{Cost Function and Constraints}
From \bequ{eqn::control}, the cost in \bequ{eqn::cost} becomes
\begin{align}
    \label{eqn::obj}
    J &= \sum_{k = 0}^{N-1} \bm{v}_k^T \bm{v}_k \, ,
\end{align}
which has a gradient with respect to any given $\bm{v}_{k}$ as
\begin{align}
    \label{eqn::objderiv}
    \nabla_{\bm{v}_{k}} J = 2\bm{v}_k \, .
\end{align}

The terminal constraint of $\bm{x}_N \sim p(\bm{x}_N)$ may be mechanized in various ways. Since the final node does not have a control update, $\bm{x}_N$ is modeled in the optimization by a GM expressed as
\begin{align}
    \label{eqn::xN}
    p(\bm{x}_N) = \sum_{\ell = 1}^{L}w_{x,N}^{\l} p_g(\bm{x}_N; \bm{m}_{x,N}^{\l}, \bm{P}_{xx,N}^{\l} ) \, .
\end{align}
Constraints on the mean and covariance of the GM do not theoretically guarantee a constrained terminal pdf as they do for a truly Gaussian problem. Furthermore, the mean of the GM does not necessarily coincide with the mode (i.e., maximum likelihood); there is, however, no closed-form solution for the mode of a GM. The general objective is to constrain the spread of the final pdf to a desired region of the state space. Therefore, this work relaxes the original problem by considering terminal constraints on the mean and covariance of the final GM pdf. The terminal constraint on the mean is given by
\begin{align}
    \label{eqn::conmean}
    \bm{m}_{x,N} = \bm{x}_f \, ,
\end{align}
where $\bm{x}_f$ is a prescribed final state. The mean $\bm{m}_{x,N}$ is computed by \bequ{eqn::GMmean}. For the covariance, multiple constraint parameterizations are available (trace, norm, etc.). Let $\bm{\lambda}(\bm{A})$ be the array of eigenvalues of matrix $\bm{A}$. The covariance constraint is given by
\begin{align}
    \label{eqn::concov}
    \bm{\lambda}(\bm{P}_{xx,f} -\bm{P}_{xx,N}) \geq \bm{0}_{n \times 1} \, ,
\end{align}
where $\bm{P}_{xx,N}$ is computed by \bequ{eqn::GMcov} and $\bm{P}_{xx,f}$ is a prescribed final covariance which, by definition, is positive semidefinite. The constraint enforces that the covariance $\bm{P}_{xx,N}$ be smaller in an eigenvalue sense than $\bm{P}_{xx,f}$. Although these two constraints do not encapsulate all probability mass for the GM, in many well-behaved systems, the mean and covariance parameterize a large portion of the GM pdf. 

The control magnitude constraint given by \bequ{eqn::prob-u} is converted into the following deterministic constraint for each node
\begin{align}
    \label{eqn::concontrol}
    \| \bm{v}_k \| + \gamma \sqrt{\lambda_{\max}(\bm{P}_{uu,k})}\left ( \sqrt{2\log\frac{1}{\beta}} + \sqrt{m}  \right ) \leq \rho_u \, ,
\end{align}
where $m$ is the dimension of the control vector, $\gamma$ is a scaling parameter where $\gamma \geq 1$, $\beta$ is the prescribed maximum probability of constraint violation where $0 < \beta \leq 1$, and $\bm{P}_{uu,k}$ is
\begin{align}
    \label{eqn::controlcov}
    \bm{P}_{uu,k} &= \exv{[\bm{u}_k - \bm{m}_{u,k}][\bm{u}_k - \bm{m}_{u,k}]^T} \nonumber \\
    &= \exv{(\bm{G}_k [\bm{x}_{k_-} - \bm{m}_{x,k_-}] + \delta \bm{u}_k)(\bm{G}_k [\bm{x}_{k_-} - \bm{m}_{x,k_-}] + \delta \bm{u}_k)^T} \nonumber \\
    &= \exv{ \bm{G}_k [\bm{x}_{k_-} - \bm{m}_{x,k_-}][\bm{x}_{k_-} - \bm{m}_{x,k_-}]^T \bm{G}_k^T + \delta \bm{u}_k\delta \bm{u}_k^T }  \nonumber \\
    &= \bm{G}_k\exv{[\bm{x}_{k_-} - \bm{m}_{x,k_-}][\bm{x}_{k_-} - \bm{m}_{x,k_-}]^T}\bm{G}_k^T + \exv{\delta \bm{u}_k\delta \bm{u}_k^T} \nonumber \\
    &= \bm{G}_k \bm{P}_{xx,k} \bm{G}_k^T + \bm{P}_{\delta u \delta u, k} \, ,
\end{align}
where it has been assumed that the state and control execution error are uncorrelated. The conversion in \bequ{eqn::concontrol} was originally proposed (for $\gamma = 1$) by Ridderhof \cite{ridderhof2019nonlinear} and proven to be a strong bound of the probability evaluation assuming that the control $\bm{u}_k$ is Gaussian. For the GM case, the theoretical guarantees do not hold. However, the parameter $\gamma$, introduced as a novelty in this work, allows for the designer to conveniently tune the constraint function in order to meet requirements. As $\gamma$ grows, the control distribution at each node must be more compact to satisfy the constraint. 

\subsubsection*{Summary of the Nonlinear Program}
From the decomposition and relaxations made in the previous section, the final nonlinear program is summarized as
\begin{subequations}
    \label{eqn::NLP}
    \begin{align}
        \label{eqn::NLP-a}
        &\min_{\bm{v}_{0:N-1}, \bm{G}_{0:N-1}} \qquad \qquad \qquad \sum_{k = 0}^{N-1} \bm{v}_k^T\bm{v}_k  \\
        &\text{subject to node constraints,} \nonumber \\
        \label{eqn::NLP-b}
        & \hspace{5mm} \| \bm{v}_k \| + \gamma \sqrt{\lambda_{\max}(\bm{P}_{uu,k})}\left ( \sqrt{2\log\frac{1}{\beta}} + \sqrt{m}  \right ) \leq \rho_u \qquad \forall k \in [0, N-1] \\
        &\text{and terminal constraints,} \nonumber \\
        \label{eqn::NLP-c}
        & \hspace{4cm} \bm{m}_{x,N} = \bm{x}_f \\
        \label{eqn::NLP-d}
        & \hspace{20.5mm} \bm{\lambda}(\bm{P}_{xx,f} -\bm{P}_{xx,N}) \geq \bm{0}_{n \times 1} \, .
    \end{align}
\end{subequations}
The optimization variables are the set of nominal controls $\bm{v}_{0:N-1}$ and the set of feedback gains $\bm{G}_{0:N-1}$. Although the cost in \bequ{eqn::NLP-a} is convex in the optimization variables, the program in its current form is nonconvex due to the constraints in  \bequt{eqn::NLP-b}{eqn::NLP-d}. It has been shown by various authors\cite{ridderhof2020chance, skaf2010design, bakolas2016optimal} that, when the dynamics are linear, the problem posed by \bequg{eqn::NLP} can be transformed into a convex program. For the current analysis, the objective is to solve a nonlinear program obtained directly from \bequg{eqn::NLP}, which can be solved using a variety of nonlinear programming methods such as sequential quadratic programming \cite{boggs1995sequential}. 

\section*{Earth-To-Mars Transfer Problem}
To evaluate the proposed methodology, a three-dimensional, finite-burn Earth-to-Mars transfer problem is considered and investigated for multiple scenarios. This problem has been used as a test case for previous works in the mission design literature \cite{ridderhof2020chance, marmo2022chance, ozaki2020tube}. The spacecraft is simulated as a point mass with constant mass subjected to two-body dynamics with an additive perturbing acceleration modeled as white-noise, where the primary body is the Sun. The control is assumed to be a series of impulses added onto the spacecraft velocity vector to approximate a finite thrust. Initial conditions are the spacecraft position $\bm{r}_0$ and velocity $\bm{v}_0$ in heliocentric coordinates, and the initial Gaussian distribution of the spacecraft state given by mean $\bm{m}_{x,0}$ and covariance $\bm{P}_{xx,0}$. The problem is nondimensionalized in distance by $\text{DU} = 1$ AU (approximately 150 million km), in time by $\text{TU} = 58$ days, and in velocity by $\text{VU} = \text{DU} / \text{TU}$. The initial mean is set to the initial state given in \tbl{tbl::simparams}, and the initial and final covariances are
\begin{align*}
    \bm{P}_{xx,0} &= \begin{bmatrix}
        10^{-4}\bm{I}_3 \hspace{1mm} \text{DU}^2 & \bm{0}_{3\times3} \\
        \bm{0}_{3\times3} &  10^{-5}\bm{I}_3 \hspace{1mm} \text{VU}^2
    \end{bmatrix} \\
    \bm{P}_{xx,f} &= \begin{bmatrix}
        10^{-6}\bm{I}_3 \hspace{1mm} \text{DU}^2 & \bm{0}_{3\times3} \\
        \bm{0}_{3\times3} &  10^{-7}\bm{I}_3 \hspace{1mm} \text{VU}^2
    \end{bmatrix} \, .
\end{align*}
Three different scenarios are considered in this analysis, with the difference between them being the total number of components in the GM during propagation. The number of components for each scenario (L), and how they are generated, is shown in \tbl{tbl::scenarios}. Otherwise, the problem data remains identical between each scenario and is given in \tbl{tbl::simparams}. It is important to note that the problem data in \tbl{tbl::simparams} is identical to the test data from previous works\cite{marmo2022chance, ridderhof2020chance} except for the number of nodes $N$. The first scenario runs the optimization with a single Gaussian. Scenario two uses a GM that has components generated from splitting in the velocity dimensions. Scenario three uses a GM whose components originate from splitting in the position and velocity dimensions. The single Gaussian scenario is included to compare and contrast the performance alongside the GM scenarios.
\begin{table}[htbp]
    \centering
    \caption{Simulation parameters for the Earth-to-Mars transfer (Matlab notation is used for arrays).}
    \label{tbl::simparams}
    \begin{NiceTabular}{ccc} \hline \hline
        Quantity & Value & Unit \\
        \hline 
        $\mu$ & $1$ & $\frac{\text{DU}^3}{\text{TU}^2}$ \\
        $\bm{r}_0$ & [$-0.940$; $-0.345$; $0.000$] & DU \\
        $\bm{v}_0$ & [$0.328$; $-0.942$; $0.000$] & VU \\
        $\bm{r}_f$ & [$-1.154$; $1.183$; $0.053$] & DU \\
        $\bm{v}_f$ & [$-0.551$; $-0.498$; $0.003$] & VU \\
        $t_f - t_0$ & $358$ & days \\
        $N$ & $31$ & - \\
        $\rho_u$ & $0.76$ & $\frac{\text{km}}{\text{s}}$ \\
        $\beta$ & $0.05$ & - \\
        $\bm{F}_w$ & $\bm{I}_6$ & - \\
        $\bm{F}_u$ & [$\bm{0}_{3\times3}$ \hspace{0.2mm}; \hspace{0.2mm}$\bm{I}_3$]& - \\
        $\bm{P}_{\delta u \delta u}$ & $\bm{0}_{3\times3}$ & VU \\
        \hline \hline
    \end{NiceTabular}
\end{table}
\begin{table}[htbp]
    \centering
    \caption{Total components for each scenario with the splitting library used.}
    \label{tbl::scenarios}
    \begin{NiceTabular}{ccc} \hline \hline
        $L$ & Dimensions & Library \\
        \hline 
        1 & - & - \\
        27 & velocity & 3 \\
        243 & position \& velocity & 3 \\
        \hline \hline
    \end{NiceTabular}
\end{table}

For each scenario, the analysis is comprised of three steps. The first step is to run the optimization problem through a nonlinear programming solver. This work uses SNOPT\cite{gill2005snopt}, which forms iterative quadratic programs using finite differences for constraint derivatives. The optimizer is given the cost function gradient from \bequ{eqn::objderiv}. The result of this step is the optimized nominal controls $\bm{v}_{0:N-1}$ and closed-loop gains $\bm{G}_{0:N-1}$. The second step is to propagate the initial pdf through time using these optimized control parameters. The result of this step is a set of ``solved for" means $\bm{m}_{x,k}$ and covariances $\bm{P}_{xx,k}$ at each node $k$. The ``solved for" means are a by-product of the optimization and are used in \bequ{eqn::control}. The final step is a Monte Carlo simulation that initializes 50,000 samples of the initial pdf and propagates each sample separately using the dynamics and feedback law in \bequ{eqn::control}. The result of the Monte Carlo is a trajectory and control history for each sample. The optimization steps for each scenario are seeded with an initial guess for the nominal control given by equally distributing the $\Delta V$ computed by a simple Hohmann transfer. The initial guess for the closed-loop gain is set to zero. Depending on the number of GM components, computational time for the optimization steps ranged from 15 seconds ($L = 27$) to 8 minutes ($L = 243$). 

\section*{Results and Discussion}
Since each scenario is successfully optimized with all constraints satisfied, and for the sake of brevity, it suffices to show the optimization results for the single Gaussian ($L = 1$) scenario and the $L = 27$ scenario. The optimal trajectories, $\Delta V$ vectors, and covariances in the planar view along with a zoomed portion of the final conditions are given for both scenarios in \bfigs{fig::1optTraj}{fig::27optTrajzoom}. The orbits of Earth and Mars are plotted in blue and red, respectively, with the optimized trajectory (generated via the optimal means) shown in yellow. The black ellipses in \bfigt{fig::1optTraj}{fig::27optTraj} portray the 99.75\% interval covariance (assuming a 2-D Gaussian). In \bfigt{fig::1optTrajzoom}{fig::27optTrajzoom} the green dashed ellipses are the desired final covariance, and the black solid ellipses are the optimized final covariance. It is clear that both cases satisfy the terminal covariance constraint, with the $L = 27$ scenario having a slightly tigher covariance. The optimized, nominal control magnitudes across all nodes for scenarios $L = 1$ and $L = 27$ are shown in \bfigt{fig::1dvmag}{fig::27dvmag}. Both scenarios (and all scenarios) converge to similar control solutions, and the magnitude constraint is never violated by the nominal control. 
\begin{figure}[h!]
\centering
    \begin{minipage}[b]{0.47\textwidth}
        \centering
        \setlength\figurewidth{6.7cm}
        \setlength\figureheight{6.7cm}
        \input{posTraj_1}
        \caption{$\bm{L=1}$ optimized means (yellow), covariance ellipses (black), and $\bm{\Delta V}$ vectors (purple).}
        \label{fig::1optTraj}
    \end{minipage}\hspace{3mm}
    \begin{minipage}[b]{0.47\textwidth}
        \centering
        \setlength\figurewidth{6.7cm}
        \setlength\figureheight{6.7cm}
\definecolor{mycolor1}{rgb}{0.00000,0.44700,0.74100}%
\definecolor{mycolor2}{rgb}{0.85000,0.32500,0.09800}%
\definecolor{mycolor3}{rgb}{0.92900,0.69400,0.12500}%
\definecolor{mycolor4}{rgb}{0.49400,0.18400,0.55600}%
\begin{tikzpicture}

\begin{axis}[%
width=0.889\figurewidth,
height=0.889\figureheight,
at={(0\figurewidth,0\figureheight)},
scale only axis,
xmin=-1.159,
xmax=-1.15,
xlabel style={font=\color{white!15!black}},
xlabel={x-position [DU]},
ymin=1.178,
ymax=1.188,
axis background/.style={fill=white},
ylabel style={font=\color{white!15!black}, font=\footnotesize},
ylabel={$y$-position [DU]},
xlabel style={font=\color{white!15!black}, font=\footnotesize},
yticklabel style={/pgf/number format/.cd,fixed zerofill,precision=3},
xticklabel style={/pgf/number format/.cd,fixed zerofill,precision=3}
]
\addplot [color=mycolor3, mark size=2pt, mark=*, mark options={mycolor3, fill=mycolor3}, forget plot]
  table[row sep=crcr]{%
-1.15429138098565	1.18288231794485\\
};
\addplot [color=black, forget plot]
  table[row sep=crcr]{%
-1.15162233295723	1.18336043284089\\
-1.15162770661442	1.18357020000188\\
-1.15164380594818	1.18377719730461\\
-1.15167056613217	1.183980591244\\
-1.15170787941258	1.18417956282444\\
-1.15175559554198	1.18437331085761\\
-1.15181352238438	1.18456105518855\\
-1.15188142668883	1.18474203983711\\
-1.15195903502865	1.18491553604199\\
-1.15204603490245	1.18508084519521\\
-1.15214207599244	1.18523730165514\\
-1.15224677157506	1.18538427542685\\
-1.15235970007815	1.18552117469883\\
-1.15248040677847	1.18564744822603\\
-1.15260840563275	1.18576258754953\\
-1.15274318123479	1.1858661290439\\
-1.15288419089079	1.18595765578409\\
-1.15303086680465	1.18603679922422\\
-1.15318261836422	1.18610324068159\\
-1.15333883451954	1.18615671261991\\
-1.15349888624329	1.18619699972655\\
-1.15366212906365	1.18622393977955\\
-1.15382790565941	1.18623742430083\\
-1.15399554850673	1.18623739899298\\
-1.15416438256704	1.1862238639579\\
-1.15433372800515	1.1861968736964\\
-1.15450290292676	1.18615653688874\\
-1.15467122612416	1.18610301595702\\
-1.15483801981926	1.18603652641112\\
-1.15500261239271	1.18595733598099\\
-1.15516434108831	1.18586576353853\\
-1.15532255468168	1.18576217781365\\
-1.15547661610252	1.18564699590952\\
-1.15562590499987	1.18552068162299\\
-1.15576982024002	1.18538374357713\\
-1.15590778232711	1.18523673317311\\
-1.15603923573649	1.18508024236994\\
-1.15616365115168	1.18491490130086\\
-1.15628052759572	1.18474137573598\\
-1.15638939444842	1.18456036440153\\
-1.15648981334138	1.18437259616625\\
-1.15658137992317	1.18417882710656\\
-1.15666372548748	1.18397983746207\\
-1.15673651845779	1.18377642849384\\
-1.15679946572249	1.18356941925801\\
-1.15685231381517	1.18335964330769\\
-1.15689484993521	1.18314794533661\\
-1.15692690280467	1.18293517777777\\
-1.15694834335795	1.18272219737105\\
-1.15695908526152	1.18250986171335\\
-1.15695908526152	1.18229902580545\\
-1.15694834335795	1.18209053860913\\
-1.15692690280467	1.18188523962873\\
-1.15689484993521	1.1816839555308\\
-1.15685231381517	1.18148749681536\\
-1.15679946572249	1.18129665455231\\
-1.15673651845779	1.1811121971961\\
-1.15666372548748	1.18093486749139\\
-1.15658137992317	1.18076537948231\\
-1.15648981334138	1.18060441563726\\
-1.15638939444842	1.18045262410083\\
-1.15628052759572	1.18031061608397\\
-1.15616365115168	1.18017896340284\\
-1.15603923573649	1.18005819617631\\
-1.15590778232711	1.17994880069139\\
-1.15576982024002	1.17985121744508\\
-1.15562590499987	1.17976583937068\\
-1.15547661610252	1.17969301025555\\
-1.15532255468168	1.17963302335684\\
-1.15516434108831	1.17958612022062\\
-1.15500261239271	1.17955248970926\\
-1.15483801981926	1.17953226724097\\
-1.15467122612416	1.17952553424448\\
-1.15450290292676	1.17953231783119\\
-1.15433372800515	1.179552590686\\
-1.15416438256704	1.17958627117728\\
-1.15399554850673	1.17963322368557\\
-1.15382790565941	1.1796932591497\\
-1.15366212906365	1.17976613582803\\
-1.15349888624329	1.17985156027191\\
-1.15333883451954	1.17994918850726\\
-1.15318261836422	1.18005862741962\\
-1.15303086680465	1.18017943633712\\
-1.15288419089079	1.18031112880489\\
-1.15274318123479	1.18045317454384\\
-1.15260840563275	1.18060500158591\\
-1.15248040677847	1.1807659985772\\
-1.15235970007815	1.18093551723964\\
-1.15224677157506	1.18111287498141\\
-1.15214207599244	1.18129735764548\\
-1.15204603490245	1.18148822238527\\
-1.15195903502865	1.18168470065585\\
-1.15188142668883	1.18188600130855\\
-1.15181352238438	1.18209131377671\\
-1.15175559554198	1.18229981133947\\
-1.15170787941258	1.18251065445074\\
-1.15167056613217	1.18272299411972\\
-1.15164380594818	1.18293597532951\\
-1.15162770661442	1.18314874047996\\
-1.15162233295723	1.18336043284089\\
};
\addplot [color=green, dashed, forget plot]
  table[row sep=crcr]{%
-1.15083038393361	1.18288415106952\\
-1.1508373524765	1.1831036748066\\
-1.15085823004533	1.18332231459908\\
-1.1508929325735	1.18353919006167\\
-1.15094132032615	1.18375342791343\\
-1.15100319846291	1.18396416549414\\
-1.15107831782234	1.18417055423795\\
-1.15116637592531	1.18437176309023\\
-1.15126701819291	1.18456698185395\\
-1.15137983937424	1.18475542445204\\
-1.15150438517824	1.18493633209267\\
-1.15164015410291	1.18510897632463\\
-1.1517865994547	1.18527266197052\\
-1.15194313154988	1.18542672992603\\
-1.15210912008894	1.18557055981387\\
-1.15228389669461	1.18570357248184\\
-1.1524667576032	1.18582523233489\\
-1.15265696649837	1.18593504949171\\
-1.15285375747602	1.18603258175739\\
-1.15305633812835	1.18611743640391\\
-1.15326389273457	1.18618927175155\\
-1.15347558554554	1.18624779854472\\
-1.15369056414901	1.18629278111669\\
-1.15390796290204	1.18632403833852\\
-1.15412690641657	1.1863414443484\\
-1.15434651308434	1.18634492905848\\
-1.15456589862681	1.18633447843706\\
-1.15478417965583	1.18631013456509\\
-1.15500047723074	1.18627199546677\\
-1.15521392039755	1.18622021471479\\
-1.15542364969596	1.18615500081194\\
-1.1556288206201	1.18607661635161\\
-1.15582860701907	1.18598537696035\\
-1.15602220442357	1.18588165002698\\
-1.15620883328517	1.18576585322326\\
-1.15638774211532	1.18563845282205\\
-1.1565582105113	1.18549996181979\\
-1.15671955205704	1.18535093787088\\
-1.15687111708709	1.18519198104215\\
-1.15701229530257	1.18502373139666\\
-1.15714251822863	1.18484686641634\\
-1.1572612615035	1.18466209827405\\
-1.15736804698994	1.18447017096588\\
-1.15746244470047	1.18427185731537\\
-1.15754407452883	1.18406795586158\\
-1.15761260778051	1.18385928764372\\
-1.1576677684963	1.18364669289504\\
-1.15770933456346	1.18343102765956\\
-1.15773713861011	1.18321316034507\\
-1.15775106867918	1.18299396822635\\
-1.15775106867918	1.18277433391269\\
-1.15773713861011	1.18255514179397\\
-1.15770933456346	1.18233727447948\\
-1.1576677684963	1.182121609244\\
-1.15761260778051	1.18190901449532\\
-1.15754407452883	1.18170034627745\\
-1.15746244470047	1.18149644482367\\
-1.15736804698994	1.18129813117315\\
-1.1572612615035	1.18110620386498\\
-1.15714251822863	1.18092143572269\\
-1.15701229530257	1.18074457074238\\
-1.15687111708709	1.18057632109688\\
-1.15671955205704	1.18041736426816\\
-1.1565582105113	1.18026834031924\\
-1.15638774211532	1.18012984931699\\
-1.15620883328517	1.18000244891577\\
-1.15602220442357	1.17988665211205\\
-1.15582860701907	1.17978292517869\\
-1.1556288206201	1.17969168578743\\
-1.15542364969596	1.1796133013271\\
-1.15521392039755	1.17954808742425\\
-1.15500047723074	1.17949630667226\\
-1.15478417965583	1.17945816757394\\
-1.15456589862681	1.17943382370198\\
-1.15434651308434	1.17942337308056\\
-1.15412690641657	1.17942685779064\\
-1.15390796290204	1.17944426380052\\
-1.15369056414901	1.17947552102234\\
-1.15347558554554	1.17952050359432\\
-1.15326389273457	1.17957903038749\\
-1.15305633812835	1.17965086573513\\
-1.15285375747602	1.17973572038165\\
-1.15265696649837	1.17983325264732\\
-1.1524667576032	1.17994306980415\\
-1.15228389669461	1.18006472965719\\
-1.15210912008894	1.18019774232517\\
-1.15194313154988	1.18034157221301\\
-1.1517865994547	1.18049564016851\\
-1.15164015410291	1.18065932581441\\
-1.15150438517824	1.18083197004637\\
-1.15137983937424	1.181012877687\\
-1.15126701819291	1.18120132028509\\
-1.15116637592531	1.18139653904881\\
-1.15107831782234	1.18159774790109\\
-1.15100319846291	1.1818041366449\\
-1.15094132032615	1.18201487422561\\
-1.1508929325735	1.18222911207737\\
-1.15085823004533	1.18244598753996\\
-1.1508373524765	1.18266462733243\\
-1.15083038393361	1.18288415106952\\
};
\end{axis}
\end{tikzpicture}%
        \caption{$\bm{L=1}$ terminal desired covariance (green), optimized covariance (black), and mean (yellow).}
        \label{fig::1optTrajzoom}
    \end{minipage}
\end{figure}
\begin{figure}[h!]
\centering
    \begin{minipage}{0.47\textwidth}
        \centering
        \setlength\figurewidth{6.7cm}
        \setlength\figureheight{6.7cm}
        \input{posTraj_27}
        \caption{$\bm{L=27}$ optimized means (yellow), covariance ellipses (black), and $\bm{\Delta V}$ vectors (purple).}
        \label{fig::27optTraj}
    \end{minipage}\hspace{3mm}
    \begin{minipage}{0.47\textwidth}
        \centering
        \setlength\figurewidth{6.7cm}
        \setlength\figureheight{6.7cm}
\definecolor{mycolor1}{rgb}{0.00000,0.44700,0.74100}%
\definecolor{mycolor2}{rgb}{0.85000,0.32500,0.09800}%
\definecolor{mycolor3}{rgb}{0.92900,0.69400,0.12500}%
\definecolor{mycolor4}{rgb}{0.49400,0.18400,0.55600}%
\begin{tikzpicture}

\begin{axis}[%
width=0.889\figurewidth,
height=0.889\figureheight,
at={(0\figurewidth,0\figureheight)},
scale only axis,
xmin=-1.159,
xmax=-1.15,
xlabel style={font=\color{white!15!black}},
xlabel={x-position [DU]},
ymin=1.178,
ymax=1.188,
axis background/.style={fill=white},
ylabel style={font=\color{white!15!black}, font=\footnotesize},
ylabel={$y$-position [DU]},
xlabel style={font=\color{white!15!black}, font=\footnotesize},
yticklabel style={/pgf/number format/.cd,fixed zerofill,precision=3},
xticklabel style={/pgf/number format/.cd,fixed zerofill,precision=3}
]
\addplot [color=mycolor3, mark size=2pt, mark=*, mark options={mycolor3, fill=mycolor3}, forget plot]
  table[row sep=crcr]{%
-1.1542915881365	1.18288411517257\\
};
\addplot [color=black, forget plot]
  table[row sep=crcr]{%
-1.15174550441691	1.18222100143886\\
-1.15175063050716	1.18240225943735\\
-1.15176598813698	1.1825854576988\\
-1.1517915154666	1.1827698585485\\
-1.15182710970647	1.18295471946929\\
-1.15187262753116	1.1831392960915\\
-1.15192788565646	1.18332284519022\\
-1.15199266157738	1.183504627678\\
-1.15206669446417	1.18368391158093\\
-1.15214968621251	1.183859974986\\
-1.15224130264393	1.18403210894803\\
-1.15234117485139	1.1841996203443\\
-1.15244890068476	1.18436183466555\\
-1.15256404637014	1.18451809873195\\
-1.15268614825649	1.18466778332326\\
-1.1528147146826	1.18481028571246\\
-1.15294922795686	1.18494503209273\\
-1.15308914644179	1.18507147988795\\
-1.15323390673502	1.18518911993747\\
-1.15338292593796	1.18529747854635\\
-1.15353560400286	1.1853961193927\\
-1.15369132614905	1.18548464528465\\
-1.1538494653384	1.18556269975968\\
-1.15400938480022	1.18562996851997\\
-1.15417044059527	1.18568618069796\\
-1.15433198420871	1.18573110994705\\
-1.15449336516143	1.18576457535301\\
-1.15465393362928	1.18578644216247\\
-1.15481304305972	1.1857966223255\\
-1.15497005277523	1.18579507485017\\
-1.1551243305531	1.18578180596763\\
-1.15527525517117	1.18575686910698\\
-1.15542221890927	1.18572036468015\\
-1.15556462999629	1.1856724396776\\
-1.15570191499306	1.18561328707638\\
-1.15583352110135	1.18554314506315\\
-1.1559589183898	1.18546229607504\\
-1.15607760192779	1.18537106566239\\
-1.15618909381861	1.18526982117787\\
-1.15629294512376	1.18515897029729\\
-1.15638873767068	1.18503895937801\\
-1.15647608573661	1.18491027166166\\
-1.15655463760171	1.18477342532825\\
-1.15662407696535	1.18462897140966\\
-1.1566841242197	1.18447749157084\\
-1.15673453757568	1.18431959576763\\
-1.15677511403649	1.18415591979072\\
-1.15680569021505	1.18398712270548\\
-1.15682614299187	1.18381388419821\\
-1.15683639001085	1.18363690183922\\
-1.15683639001085	1.18345688827402\\
-1.15682614299187	1.1832745683537\\
-1.15680569021505	1.18309067621623\\
-1.15677511403649	1.18290595233033\\
-1.15673453757568	1.18272114051392\\
-1.1566841242197	1.18253698493893\\
-1.15662407696535	1.18235422713487\\
-1.15655463760171	1.18217360300292\\
-1.15647608573661	1.18199583985271\\
-1.15638873767068	1.1818216534737\\
-1.15629294512376	1.18165174525296\\
-1.15618909381861	1.18148679935093\\
-1.15607760192779	1.18132747994654\\
-1.1559589183898	1.1811744285628\\
-1.15583352110135	1.18102826148361\\
-1.15570191499306	1.18088956727222\\
-1.15556462999629	1.18075890440126\\
-1.15542221890927	1.180636799004\\
-1.15527525517117	1.18052374275577\\
-1.1551243305531	1.18042019089418\\
-1.15497005277523	1.18032656038604\\
-1.15481304305972	1.18024322824835\\
-1.15465393362928	1.18017053003021\\
-1.15449336516143	1.18010875846171\\
-1.15433198420871	1.18005816227513\\
-1.15417044059527	1.18001894520347\\
-1.15400938480022	1.17999126516003\\
-1.1538494653384	1.17997523360258\\
-1.15369132614905	1.17997091508454\\
-1.15353560400286	1.17997832699505\\
-1.15338292593796	1.17999743948898\\
-1.15323390673502	1.18002817560705\\
-1.15308914644179	1.18007041158577\\
-1.15294922795686	1.18012397735575\\
-1.1528147146826	1.18018865722655\\
-1.15268614825649	1.18026419075515\\
-1.15256404637014	1.18035027379469\\
-1.15244890068476	1.18044655971916\\
-1.15234117485139	1.18055266081912\\
-1.15224130264393	1.18066814986288\\
-1.15214968621251	1.18079256181684\\
-1.15206669446417	1.18092539571796\\
-1.15199266157738	1.18106611669101\\
-1.15192788565646	1.18121415810229\\
-1.15187262753116	1.18136892384129\\
-1.15182710970647	1.18152979072101\\
-1.1517915154666	1.18169611098729\\
-1.15176598813698	1.18186721492713\\
-1.15175063050716	1.18204241356535\\
-1.15174550441691	1.18222100143886\\
};
\addplot [color=green, dashed, forget plot]
  table[row sep=crcr]{%
-1.15083038393361	1.18288415106952\\
-1.1508373524765	1.1831036748066\\
-1.15085823004533	1.18332231459908\\
-1.1508929325735	1.18353919006167\\
-1.15094132032615	1.18375342791343\\
-1.15100319846291	1.18396416549414\\
-1.15107831782234	1.18417055423795\\
-1.15116637592531	1.18437176309023\\
-1.15126701819291	1.18456698185395\\
-1.15137983937424	1.18475542445204\\
-1.15150438517824	1.18493633209267\\
-1.15164015410291	1.18510897632463\\
-1.1517865994547	1.18527266197052\\
-1.15194313154988	1.18542672992603\\
-1.15210912008894	1.18557055981387\\
-1.15228389669461	1.18570357248184\\
-1.1524667576032	1.18582523233489\\
-1.15265696649837	1.18593504949171\\
-1.15285375747602	1.18603258175739\\
-1.15305633812835	1.18611743640391\\
-1.15326389273457	1.18618927175155\\
-1.15347558554554	1.18624779854472\\
-1.15369056414901	1.18629278111669\\
-1.15390796290204	1.18632403833852\\
-1.15412690641657	1.1863414443484\\
-1.15434651308434	1.18634492905848\\
-1.15456589862681	1.18633447843706\\
-1.15478417965583	1.18631013456509\\
-1.15500047723074	1.18627199546677\\
-1.15521392039755	1.18622021471479\\
-1.15542364969596	1.18615500081194\\
-1.1556288206201	1.18607661635161\\
-1.15582860701907	1.18598537696035\\
-1.15602220442357	1.18588165002698\\
-1.15620883328517	1.18576585322326\\
-1.15638774211532	1.18563845282205\\
-1.1565582105113	1.18549996181979\\
-1.15671955205704	1.18535093787088\\
-1.15687111708709	1.18519198104215\\
-1.15701229530257	1.18502373139666\\
-1.15714251822863	1.18484686641634\\
-1.1572612615035	1.18466209827405\\
-1.15736804698994	1.18447017096588\\
-1.15746244470047	1.18427185731537\\
-1.15754407452883	1.18406795586158\\
-1.15761260778051	1.18385928764372\\
-1.1576677684963	1.18364669289504\\
-1.15770933456346	1.18343102765956\\
-1.15773713861011	1.18321316034507\\
-1.15775106867918	1.18299396822635\\
-1.15775106867918	1.18277433391269\\
-1.15773713861011	1.18255514179397\\
-1.15770933456346	1.18233727447948\\
-1.1576677684963	1.182121609244\\
-1.15761260778051	1.18190901449532\\
-1.15754407452883	1.18170034627745\\
-1.15746244470047	1.18149644482367\\
-1.15736804698994	1.18129813117315\\
-1.1572612615035	1.18110620386498\\
-1.15714251822863	1.18092143572269\\
-1.15701229530257	1.18074457074238\\
-1.15687111708709	1.18057632109688\\
-1.15671955205704	1.18041736426816\\
-1.1565582105113	1.18026834031924\\
-1.15638774211532	1.18012984931699\\
-1.15620883328517	1.18000244891577\\
-1.15602220442357	1.17988665211205\\
-1.15582860701907	1.17978292517869\\
-1.1556288206201	1.17969168578743\\
-1.15542364969596	1.1796133013271\\
-1.15521392039755	1.17954808742425\\
-1.15500047723074	1.17949630667226\\
-1.15478417965583	1.17945816757394\\
-1.15456589862681	1.17943382370198\\
-1.15434651308434	1.17942337308056\\
-1.15412690641657	1.17942685779064\\
-1.15390796290204	1.17944426380052\\
-1.15369056414901	1.17947552102234\\
-1.15347558554554	1.17952050359432\\
-1.15326389273457	1.17957903038749\\
-1.15305633812835	1.17965086573513\\
-1.15285375747602	1.17973572038165\\
-1.15265696649837	1.17983325264732\\
-1.1524667576032	1.17994306980415\\
-1.15228389669461	1.18006472965719\\
-1.15210912008894	1.18019774232517\\
-1.15194313154988	1.18034157221301\\
-1.1517865994547	1.18049564016851\\
-1.15164015410291	1.18065932581441\\
-1.15150438517824	1.18083197004637\\
-1.15137983937424	1.181012877687\\
-1.15126701819291	1.18120132028509\\
-1.15116637592531	1.18139653904881\\
-1.15107831782234	1.18159774790109\\
-1.15100319846291	1.1818041366449\\
-1.15094132032615	1.18201487422561\\
-1.1508929325735	1.18222911207737\\
-1.15085823004533	1.18244598753996\\
-1.1508373524765	1.18266462733243\\
-1.15083038393361	1.18288415106952\\
};
\end{axis}
\end{tikzpicture}%
        \caption{$\bm{L=27}$ terminal desired covariance (green), optimized covariance (black), and mean (yellow).}
        \label{fig::27optTrajzoom}
    \end{minipage}
\end{figure}
\begin{figure}[h!]
\centering
    \begin{minipage}{0.47\textwidth}
        \centering
        \setlength\figurewidth{6.7cm}
        \setlength\figureheight{6.7cm}
\definecolor{mycolor1}{rgb}{0.00000,0.44700,0.74100}%
\definecolor{mycolor2}{rgb}{0.85000,0.32500,0.09800}%
\begin{tikzpicture}

\begin{axis}[%
width=0.889\figurewidth,
height=0.889\figureheight,
at={(0\figurewidth,0\figureheight)},
scale only axis,
xmin=0,
xmax=30,
xlabel style={font=\color{white!15!black}},
xlabel={Node},
ymin=0,
ymax=0.8,
ylabel style={font=\color{white!15!black}},
ylabel={$\Delta V$ [km/s]},
axis background/.style={fill=white},
legend style={legend cell align=left, align=left, draw=none, fill=none, anchor=north west, at={(axis cs:3.2, 0.72)}},
ylabel style={font=\color{white!15!black}, font=\footnotesize},
xlabel style={font=\color{white!15!black}, font=\footnotesize}
]
\addplot[ycomb, color=mycolor1, mark=o, mark options={solid, mycolor1}] table[row sep=crcr] {%
1	0.551088894274319\\
2	0.413976397351536\\
3	0.406844636901157\\
4	0.329944592468696\\
5	0.367691301598007\\
6	0.32878755779296\\
7	0.364065603164157\\
8	0.367786225660106\\
9	0.374848866116997\\
10	0.405086252590166\\
11	0.420532414818662\\
12	0.365913697389642\\
13	0.342427664789147\\
14	0.313698574223425\\
15	0.255781193750609\\
16	0.202396373178748\\
17	0.142714217283972\\
18	0.11800489861773\\
19	0.0744915316727594\\
20	0.0812302711648387\\
21	0.112667419659125\\
22	0.152247421826581\\
23	0.20465643705437\\
24	0.261844599079068\\
25	0.321044231529399\\
26	0.391873540728036\\
27	0.469651522506516\\
28	0.553688297263822\\
29	0.650026609641695\\
30	0.734187601311722\\
};
\addplot[forget plot, color=white!15!black] table[row sep=crcr] {%
0	0\\
30	0\\
};
\addlegendentry{\scriptsize $\Delta V$s}

\addplot [color=red, dashed]
  table[row sep=crcr]{%
0	0.76\\
30	0.76\\
};
\addlegendentry{\scriptsize Max $\Delta V$}

\addplot [color=mycolor2, mark=o, mark options={solid, mycolor2}]
  table[row sep=crcr]{%
1	0.0415536558085961\\
2	0.0727686429134996\\
3	0.103445875180402\\
4	0.128324627363216\\
5	0.156049587799114\\
6	0.180841096306409\\
7	0.208292668861783\\
8	0.236024786839496\\
9	0.264289447812894\\
10	0.294834092946945\\
11	0.326543422246247\\
12	0.354134346314862\\
13	0.379954357332734\\
14	0.403608113952331\\
15	0.422894736850839\\
16	0.438155993819287\\
17	0.448917048026014\\
18	0.45781494996495\\
19	0.463431821398843\\
20	0.469556812782796\\
21	0.478052253765186\\
22	0.48953213846266\\
23	0.504963810617058\\
24	0.524707631366339\\
25	0.548915271249692\\
26	0.578463640670838\\
27	0.613876689126328\\
28	0.65562634724645\\
29	0.704640185638254\\
30	0.76\\
};
\addlegendentry{\scriptsize Normalized Sum}

\end{axis}
\end{tikzpicture}%
        \caption{$\bm{L=1}$ nominal $\bm{\Delta V}$ magnitudes.}
        \label{fig::1dvmag}
    \end{minipage}\hspace{4mm}
    \begin{minipage}{0.47\textwidth}
        \centering
        \setlength\figurewidth{6.7cm}
        \setlength\figureheight{6.7cm}
\definecolor{mycolor1}{rgb}{0.00000,0.44700,0.74100}%
\definecolor{mycolor2}{rgb}{0.85000,0.32500,0.09800}%
\begin{tikzpicture}

\begin{axis}[%
width=0.9\figurewidth,
height=0.9\figureheight,
at={(0\figurewidth,0\figureheight)},
scale only axis,
xmin=0,
xmax=30,
xlabel style={font=\color{white!15!black}},
xlabel={Node},
ymin=0,
ymax=0.8,
axis background/.style={fill=white},
legend style={legend cell align=left, align=left, draw=none, fill=none, anchor=north west, at={(axis cs:3.2, 0.72)}},
ylabel style={font=\color{white!15!black}, font=\footnotesize},
ylabel={$\Delta V$ [km/s]},
xlabel style={font=\color{white!15!black}, font=\footnotesize}
]
\addplot[ycomb, color=mycolor1, mark=o, mark options={solid, mycolor1}] table[row sep=crcr] {%
1	0.549309605877111\\
2	0.449983437250243\\
3	0.37763303512484\\
4	0.340673253504005\\
5	0.338764402738433\\
6	0.359198435810049\\
7	0.387576887006895\\
8	0.410964514359684\\
9	0.42162987084677\\
10	0.416977443419626\\
11	0.39670121998892\\
12	0.363219909702793\\
13	0.320397016896271\\
14	0.272706894465536\\
15	0.224032292775292\\
16	0.177329552135475\\
17	0.134666855646764\\
18	0.0988829215042199\\
19	0.0770511562045969\\
20	0.0807808510268586\\
21	0.109124361259796\\
22	0.151404390728842\\
23	0.201977687661404\\
24	0.259027988976148\\
25	0.322007543003954\\
26	0.390748573275398\\
27	0.465201508423534\\
28	0.545307047166013\\
29	0.631025786648804\\
30	0.722271018183972\\
};
\addplot[forget plot, color=white!15!black] table[row sep=crcr] {%
0	0\\
30	0\\
};
\addlegendentry{\scriptsize $\Delta V$s}

\addplot [color=red, dashed]
  table[row sep=crcr]{%
0	0.76\\
30	0.76\\
};
\addlegendentry{\scriptsize Max $\Delta V$}

\addplot [color=mycolor2, mark=o, mark options={solid, mycolor2}]
  table[row sep=crcr]{%
1	0.0417618315461877\\
2	0.0759722882794407\\
3	0.104682230778949\\
4	0.130582267612291\\
5	0.156337182089742\\
6	0.183645615089017\\
7	0.213111549223799\\
8	0.244355551943805\\
9	0.276410399433734\\
10	0.308111541311355\\
11	0.338271162308499\\
12	0.365885332024335\\
13	0.390243846975404\\
14	0.410976670989957\\
15	0.428008958012953\\
16	0.441490620822432\\
17	0.451728807958072\\
18	0.459246484449566\\
19	0.465104378379379\\
20	0.471245826219809\\
21	0.479542118772786\\
22	0.491052794343214\\
23	0.506408357176919\\
24	0.526101228238398\\
25	0.550582185104354\\
26	0.580289249971687\\
27	0.615656676322797\\
28	0.657114209198694\\
29	0.70508859799761\\
30	0.76\\
};
\addlegendentry{\scriptsize Normalized Cumulative Sum}

\end{axis}
\end{tikzpicture}%
        \caption{$\bm{L=27}$ nominal $\bm{\Delta V}$ magnitudes.}
        \label{fig::27dvmag}
    \end{minipage}
\end{figure}

\subsection*{Monte Carlo Analysis}
\begin{figure}[h!]
\centering
    \begin{minipage}{0.47\textwidth}
        \centering
        \includegraphics[scale=0.9]{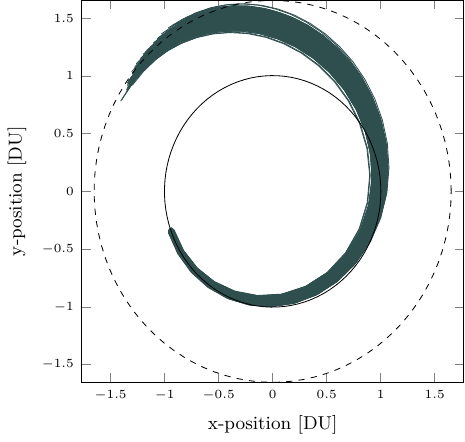}
        \caption{$\bm{L=27}$ open loop MC trajectories.}
        \label{fig::1OLtraj}
    \end{minipage}\hfill
    \begin{minipage}{0.47\textwidth}
        \centering
        \includegraphics[scale=0.9]{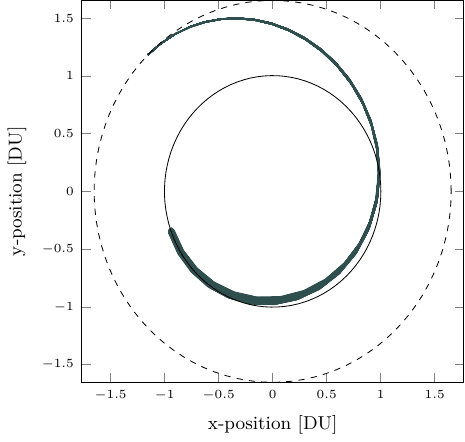}
        \caption{$\bm{L=27}$ closed-loop MC trajectories.}
        \label{fig::1CLtraj}
    \end{minipage}
\end{figure}

The objective behind the Monte Carlo (MC) analysis is to determine if the optimized closed-loop gains enable a statistically significant number of samples to reach the desired pdf while satisfying control magnitude constraints. To showcase the difference in nominal (open loop) control and closed-loop control, the $L = 27$ sample trajectories (planar) are shown in \bfigt{fig::1OLtraj}{fig::1CLtraj} as a representative example. From visual inspection it is clear the closed-loop control enables the trajectory dispersions to tighten as they reach the Mars orbit. Conversely, the nominal control does not keep the dispersions from spreading, and the terminal states are significantly far from the desired point. 

\begin{figure}[h!]
\centering
    \begin{minipage}{0.47\textwidth}
        \centering
        \includegraphics[scale=0.93]{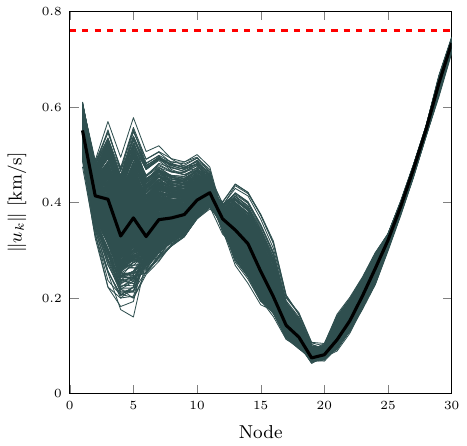}
        \caption{$\bm{L=1}$ closed-loop control magnitude histories (dark blue), nominal control magnitude (black), and maximum magnitude constraint (red).}
        \label{fig::1uCL}
    \end{minipage}\hspace{8mm}
    \begin{minipage}{0.47\textwidth}
        \centering
        \includegraphics[scale=0.93]{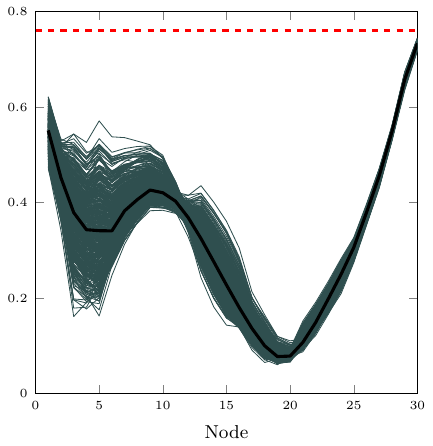}
        \caption{$\bm{L=27}$ closed-loop control magnitude histories (dark blue), nominal control magnitude (black), and maximum magnitude constraint (red).}
        \label{fig::27uCL}
    \end{minipage}
\end{figure}

The closed-loop control magnitude histories for $L = 1$, $L = 27$, and $L = 243$ are shown in \bfigs{fig::1uCL}{fig::243uCL}. Visually, the control magnitude constraint is not activated in any scenario, and all scenarios exhibit similar control histories with slight differences from one another. Since this analysis uses $\beta = 0.05$ and $\gamma = 1$ (see \bequ{eqn::concontrol}), it is expected that, if the closed-loop control magnitudes did exceed the constraint, only up to 5\% of sample histories would do so \textit{if the underlying closed-loop control distribution is truly Gaussian}. Since the continuous time dynamics are nonlinear, the state distribution becomes non-Gaussian and in turn the affine feedback law (\bequ{eqn::control}) yields a non-Gaussian control. Therefore, it could be expected that more than 5\% of sample histories to break the control constraint and the percentage could be decreased by increasing $\gamma$. An example illustrating the utility of gamma in enforcing a chance constraint, when the default unity value results in violation, is a matter for future investigation.  In the present case, gamma = 1, a value consistent with the literature, led to satisfactory results. 
\begin{figure}[h!]
    \centering
    \includegraphics[scale=1.1]{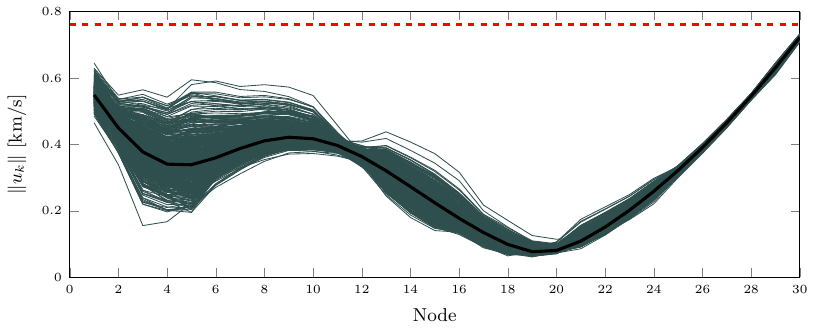}
    \caption{$\bm{L=243}$ closed-loop control magnitude histories (dark blue), nominal control magnitude (black), and maximum magnitude constraint (red).}
    \label{fig::243uCL}
\end{figure}

\begin{figure}[h!]
    \centering
    \setlength\figurewidth{14.5cm}
    \setlength\figureheight{6.3cm}
    \input{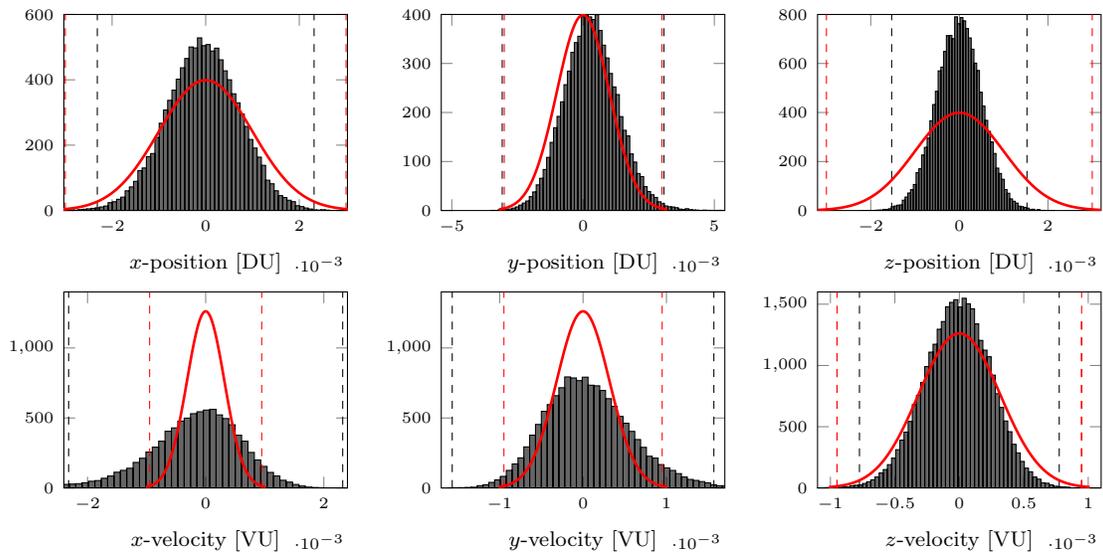}
    \caption{$\bm{L=1}$ final state error histograms (black) and desired final distributions (red).}
    \label{fig::1finalpdf}
\end{figure}

The MC performance of the terminal state and covariance constraints for $L = 1$, $L = 27$, and $L = 243$ is illustrated in the histograms of \bfigs{fig::1finalpdf}{fig::243finalpdf}. The histogram (black) for each dimension is computed by subtracting each $i$-th sample final state $\bm{x}_f^{(i)}$ from the desired final state $\bm{x}_f$. Overlaying the histograms are Gaussian distributions with zero-mean and variance given by the diagonals of the prescribed final covariance $\bm{P}_{xx,f}$. The red dashed lines indicate the prescribed $3\sigma$, and the black dashed lines the MC $3\sigma$. Every scenario is plotted with identical axes limits for comparison. Ideal performance is characterized by a significant majority of the histogram being contained within the desired distribution. To be exact, the nonlinear program as implemented requires the terminal ``solved for" covariance to be less than the prescribed covariance in the principal directions. This roughly corresponds to the black dashed lines being contained within the red dashed lines in \bfigs{fig::1finalpdf}{fig::243finalpdf}. 
\begin{figure}[h!]
    \centering
    \setlength\figurewidth{14.5cm}
    \setlength\figureheight{6.3cm}
    \input{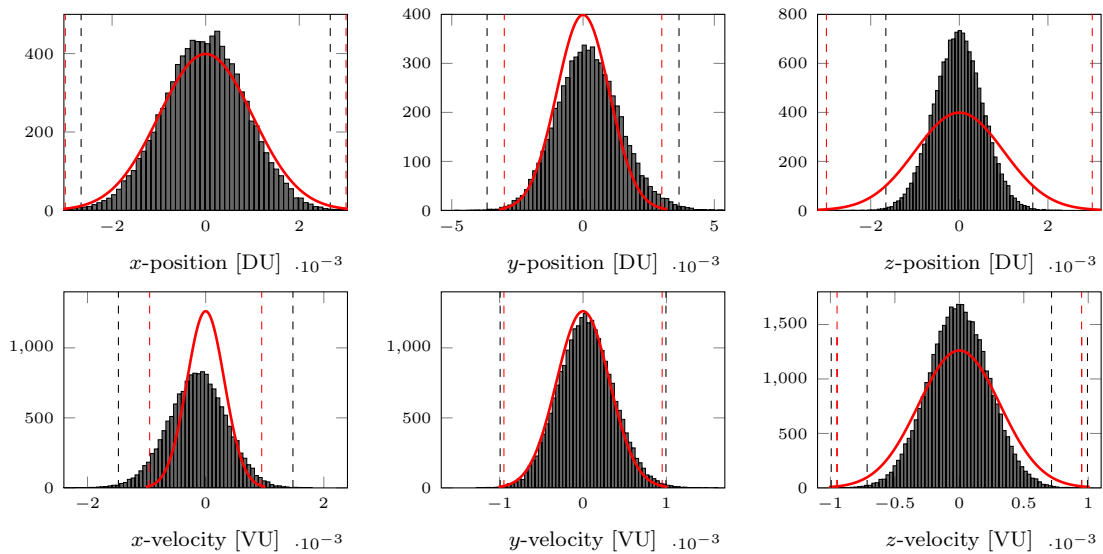}
    \caption{$\bm{L=27}$ final state error histograms (black) and desired final distributions (red).}
    \label{fig::27finalpdf}
\end{figure}
Since this is a nonlinear, non-Gaussian analysis, it is not guaranteed that the MC $3\sigma$ falls within the prescribed $3\sigma$, even though the optimization covariance constraint is satisfied. The reason is two-fold. First and foremost, the recursive splitting and merging procedure from \bfig{fig::flowchart} results in a Gaussian approximation at the end of each continuous time propagation interval. These intervals are roughly 12 days long, which is enough time to result in a significantly non-Gaussian distribution where the mean and mode do not correspond and $3\sigma$ does not encapsulate most of the information. Thus, at every node update, there is an \textit{approximation} of the true state pdf that the control acts on. A natural extension to the methodology is to remove the GM collapsing step before the node such that the control becomes distributed by a GM (due to feedback). Even then, there would not be guarantees that the MC $3\sigma$ falls within the prescribed $3\sigma$ due to the Taylor series linearization in the uncertainty propagation. This, however, serves as motivation for a GM-based approach since a higher number of components generally yields closer approximations to the true state pdf. More generally, considering that the problem is fundamentally non-Gaussian, it may be insufficient to analyze results in terms of standard deviations. Instead, an information theoretic approach may be more suitable for rigorous validation. With a practical viewpoint in mind, the objective then becomes determining how well the GM model for the optimization performs in Monte Carlo versus the single Gaussian model.
\begin{figure}[h!]
    \centering
    \setlength\figurewidth{14.5cm}
    \setlength\figureheight{6.3cm}
    \input{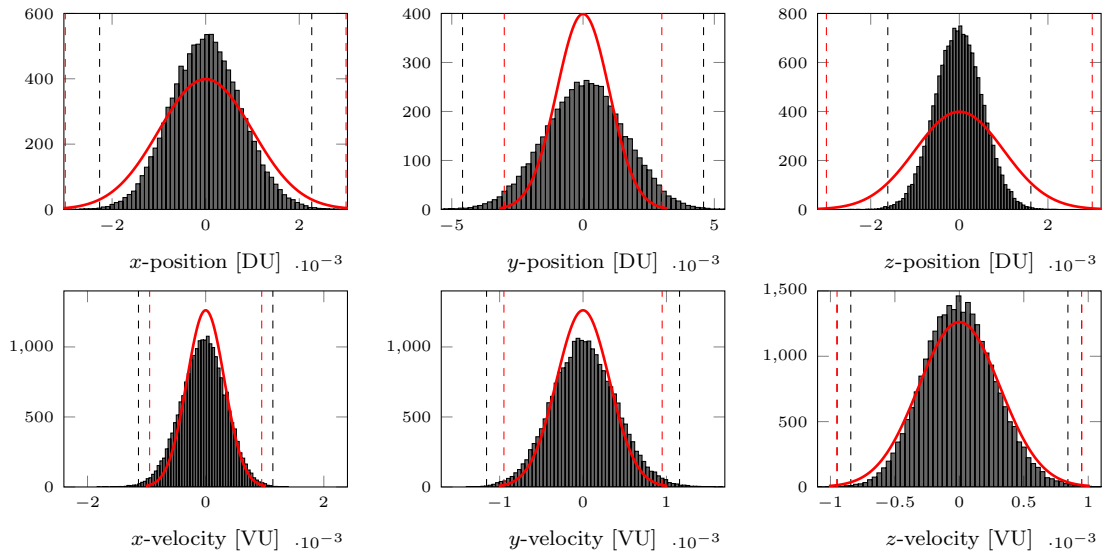}
    \caption{$\bm{L=243}$ final state error histograms (black) and desired final distributions (red).}
    \label{fig::243finalpdf}
\end{figure}

\begin{table}[h!]
\centering
    \caption{Percent of samples that fall within 99.9\% probability gate for 1D Gaussian.}
    \label{tbl::percentiles}
    \begin{NiceTabular}{ccccccc} \hline \hline
        $L$ & $x$ & $y$ & $z$ & $\dot{x}$ & $\dot{y}$ & $\dot{z}$ \\
        \hline 
        1 & 99.99 & 99.56 & 100 & 83.54 & 94.98 & 99.99 \\
        27 & 99.98 & 98.91 & 100 & 95.65 & 99.75 & 100 \\
        243 & 100 & 96.96 & 100 & 99.25 & 99.25 & 99.96 \\
        \hline \hline
    \end{NiceTabular}
\end{table}

\bfigs{fig::1finalpdf}{fig::243finalpdf} showcase non-Gaussian, heavy-tailed behavior primarily in the planar velocity states (e.g., the x-velocity histogram in \bfigt{fig::1finalpdf}{fig::27finalpdf}). Additionally, in the planar velocity histograms of \bfig{fig::1finalpdf} there is non-symmetric behavior that is not well captured by a Gaussian. Along with the figures, a quantitative comparison for terminal constraint satisfaction is given in \tbl{tbl::percentiles}. This table shows the percentage of MC samples at the final node that fall within the 99.9\% probability gate of the desired Gaussian distribution. The table shows significant improvement in overall constraint satisfaction (83\% $\rightarrow$ 99\%) when moving from the single Gaussian model to the GM model, specifically in the planar velocity directions. Overall, it is clear the GM modeling increases true feasibility by a nonnegligible amount. The underlying reasoning is because the ``solved for" means for the Gaussian mixture scenarios more accurately approximate the true state distribution mean than the single Gaussian scenario.

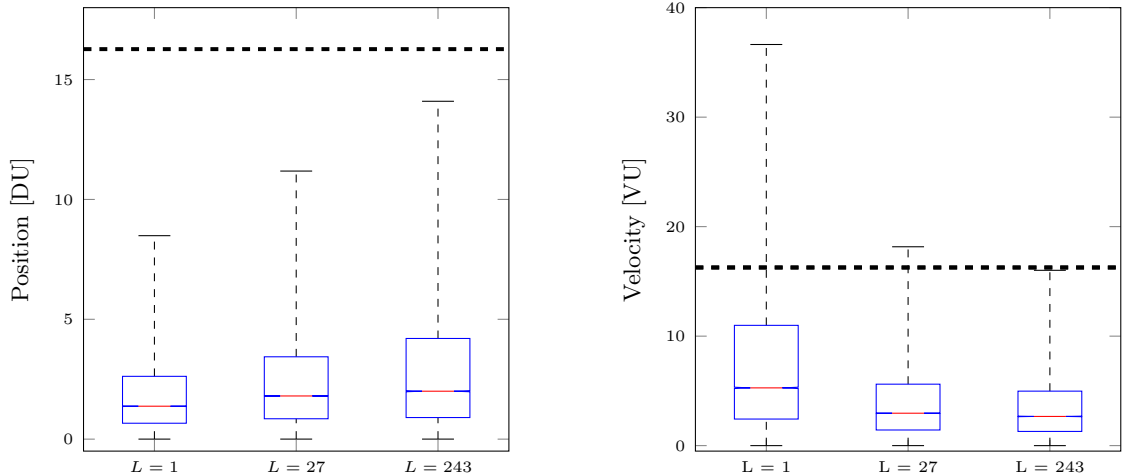
\begin{figure}[h!]
    \centering
    \setlength\figurewidth{14.5cm}
    \setlength\figureheight{5.9cm}
\begin{tikzpicture}

\begin{axis}[%
width=0.39\figurewidth,
height=\figureheight,
at={(0\figurewidth,0\figureheight)},
scale only axis,
xmin=0.5,
xmax=3.5,
xtick={1,2,3},
xticklabels={{$L = 1$},{$L = 27$},{$L = 243$}},
ymin=-0.5,
ymax=18,
ylabel style={font=\color{white!15!black}},
ylabel={Position [DU]},
axis background/.style={fill=white},
ylabel style={font=\color{white!15!black}, font=\footnotesize},
xlabel style={font=\color{white!15!black}, font=\footnotesize}
]
\addplot [color=black, dashed, forget plot]
  table[row sep=crcr]{%
1	2.61960057137693\\
1	8.48734174199385\\
};
\addplot [color=black, dashed, forget plot]
  table[row sep=crcr]{%
2	3.43489953155748\\
2	11.1828865775806\\
};
\addplot [color=black, dashed, forget plot]
  table[row sep=crcr]{%
3	4.19943484850775\\
3	14.0945824467079\\
};
\addplot [color=black, dashed, forget plot]
  table[row sep=crcr]{%
1	0.000665682858010941\\
1	0.663182553517055\\
};
\addplot [color=black, dashed, forget plot]
  table[row sep=crcr]{%
2	0.00119150360571574\\
2	0.849868906504414\\
};
\addplot [color=black, dashed, forget plot]
  table[row sep=crcr]{%
3	0.00108090936223526\\
3	0.90007373582443\\
};
\addplot [color=black, forget plot]
  table[row sep=crcr]{%
0.887499999999999	8.48734174199385\\
1.1125	8.48734174199385\\
};
\addplot [color=black, forget plot]
  table[row sep=crcr]{%
1.8875	11.1828865775806\\
2.1125	11.1828865775806\\
};
\addplot [color=black, forget plot]
  table[row sep=crcr]{%
2.8875	14.0945824467079\\
3.1125	14.0945824467079\\
};
\addplot [color=black, forget plot]
  table[row sep=crcr]{%
0.8875	0.000665682858010941\\
1.1125	0.000665682858010941\\
};
\addplot [color=black, forget plot]
  table[row sep=crcr]{%
1.8875	0.00119150360571574\\
2.1125	0.00119150360571574\\
};
\addplot [color=black, forget plot]
  table[row sep=crcr]{%
2.8875	0.00108090936223526\\
3.1125	0.00108090936223526\\
};
\addplot [color=blue, forget plot]
  table[row sep=crcr]{%
0.8875	1.37358215784082\\
0.775	1.38731866459709\\
0.775	2.61960057137693\\
1.225	2.61960057137693\\
1.225	1.38731866459709\\
1.1125	1.37358215784082\\
1.225	1.35984565108456\\
1.225	0.663182553517055\\
0.775	0.663182553517055\\
0.775	1.35984565108456\\
0.8875	1.37358215784082\\
};
\addplot [color=blue, forget plot]
  table[row sep=crcr]{%
1.8875	1.80070508850496\\
1.775	1.81885524369779\\
1.775	3.43489953155747\\
2.225	3.43489953155747\\
2.225	1.81885524369779\\
2.1125	1.80070508850496\\
2.225	1.78255493331214\\
2.225	0.849868906504414\\
1.775	0.849868906504414\\
1.775	1.78255493331214\\
1.8875	1.80070508850496\\
};
\addplot [color=blue, forget plot]
  table[row sep=crcr]{%
2.8875	1.99940007216849\\
2.775	2.02256572276156\\
2.775	4.19943484850775\\
3.225	4.19943484850775\\
3.225	2.02256572276156\\
3.1125	1.99940007216849\\
3.225	1.97623442157541\\
3.225	0.900073735824431\\
2.775	0.900073735824431\\
2.775	1.97623442157541\\
2.8875	1.99940007216849\\
};
\addplot [color=red, forget plot]
  table[row sep=crcr]{%
0.8875	1.37358215784082\\
1.1125	1.37358215784082\\
};
\addplot [color=red, forget plot]
  table[row sep=crcr]{%
1.8875	1.80070508850496\\
2.1125	1.80070508850496\\
};
\addplot [color=red, forget plot]
  table[row sep=crcr]{%
2.8875	1.99940007216848\\
3.1125	1.99940007216848\\
};
\addplot [color=black, dashed, line width=1.5pt, forget plot]
  table[row sep=crcr]{%
0.5	16.27\\
3.5	16.27\\
};
\end{axis}

\begin{axis}[%
width=0.39\figurewidth,
height=\figureheight,
at={(0.561\figurewidth,0\figureheight)},
scale only axis,
xmin=0.5,
xmax=3.5,
xtick={1,2,3},
xticklabels={{L = 1},{L = 27},{L = 243}},
ymin=-0.5,
ymax=40,
ylabel style={font=\color{white!15!black}},
ylabel={Velocity [VU]},
axis background/.style={fill=white},
ylabel style={font=\color{white!15!black}, font=\footnotesize},
xlabel style={font=\color{white!15!black}, font=\footnotesize}
]
\addplot [color=black, dashed, forget plot]
  table[row sep=crcr]{%
1	10.9855115142801\\
1	36.6344333804981\\
};
\addplot [color=black, dashed, forget plot]
  table[row sep=crcr]{%
2	5.61321863893858\\
2	18.1621203575933\\
};
\addplot [color=black, dashed, forget plot]
  table[row sep=crcr]{%
3	4.97729323840672\\
3	16.0122736537545\\
};
\addplot [color=black, dashed, forget plot]
  table[row sep=crcr]{%
1	0.000887838006783159\\
1	2.42803142088759\\
};
\addplot [color=black, dashed, forget plot]
  table[row sep=crcr]{%
2	0.00151004834450852\\
2	1.42866654221136\\
};
\addplot [color=black, dashed, forget plot]
  table[row sep=crcr]{%
3	0.000224232249869782\\
3	1.29816358594335\\
};
\addplot [color=black, forget plot]
  table[row sep=crcr]{%
0.887500000000003	36.6344333804981\\
1.1125	36.6344333804981\\
};
\addplot [color=black, forget plot]
  table[row sep=crcr]{%
1.8875	18.1621203575933\\
2.1125	18.1621203575933\\
};
\addplot [color=black, forget plot]
  table[row sep=crcr]{%
2.8875	16.0122736537545\\
3.1125	16.0122736537545\\
};
\addplot [color=black, forget plot]
  table[row sep=crcr]{%
0.8875	0.000887838006783159\\
1.1125	0.000887838006783159\\
};
\addplot [color=black, forget plot]
  table[row sep=crcr]{%
1.8875	0.00151004834450852\\
2.1125	0.00151004834450852\\
};
\addplot [color=black, forget plot]
  table[row sep=crcr]{%
2.8875	0.000224232249869782\\
3.1125	0.000224232249869782\\
};
\addplot [color=blue, forget plot]
  table[row sep=crcr]{%
0.887499999999999	5.27466301059002\\
0.775	5.33474724721349\\
0.775	10.9855115142801\\
1.225	10.9855115142801\\
1.225	5.33474724721349\\
1.1125	5.27466301059002\\
1.225	5.21457877396655\\
1.225	2.42803142088759\\
0.775	2.42803142088759\\
0.775	5.21457877396655\\
0.887499999999999	5.27466301059002\\
};
\addplot [color=blue, forget plot]
  table[row sep=crcr]{%
1.8875	2.96186541988206\\
1.775	2.99124622072519\\
1.775	5.61321863893858\\
2.225	5.61321863893858\\
2.225	2.99124622072519\\
2.1125	2.96186541988206\\
2.225	2.93248461903893\\
2.225	1.42866654221136\\
1.775	1.42866654221136\\
1.775	2.93248461903893\\
1.8875	2.96186541988206\\
};
\addplot [color=blue, forget plot]
  table[row sep=crcr]{%
2.8875	2.66767669592517\\
2.775	2.69350879768813\\
2.775	4.97729323840672\\
3.225	4.97729323840672\\
3.225	2.69350879768813\\
3.1125	2.66767669592517\\
3.225	2.64184459416221\\
3.225	1.29816358594335\\
2.775	1.29816358594335\\
2.775	2.64184459416221\\
2.8875	2.66767669592517\\
};
\addplot [color=red, forget plot]
  table[row sep=crcr]{%
0.8875	5.27466301059002\\
1.1125	5.27466301059002\\
};
\addplot [color=red, forget plot]
  table[row sep=crcr]{%
1.8875	2.96186541988206\\
2.1125	2.96186541988206\\
};
\addplot [color=red, forget plot]
  table[row sep=crcr]{%
2.8875	2.66767669592517\\
3.1125	2.66767669592517\\
};
\addplot [color=black, dashed, line width=1.5pt, forget plot]
  table[row sep=crcr]{%
0.5	16.27\\
3.5	16.27\\
};
\addplot [color=black, dashed, line width=1.5pt, forget plot]
  table[row sep=crcr]{%
0.5	16.27\\
3.5	16.27\\
};
\end{axis}
\end{tikzpicture}%
    \caption{Box-and-whisker comparison for terminal Mahalanobis distances (squared) of position and velocity MC samples. The horizontal dashed line (black) indicates the 99.9\% probability gate for a 3D Gaussian.}
    \label{fig::xfcompare}
\end{figure}
An alternative visualization using the Mahalanobis distances of the final samples is shown in \bfig{fig::xfcompare} via box-and-whisker plots. In this illustration, the 99.9\% probability gate is computed for a 3D Gaussian distribution and compared to the squared Mahalanobis distances of the samples from the desired state at the final time. The Gaussian and GM models both perform well in the position states; however, the GM models outperform the single Gaussian model significantly in the velocity states.

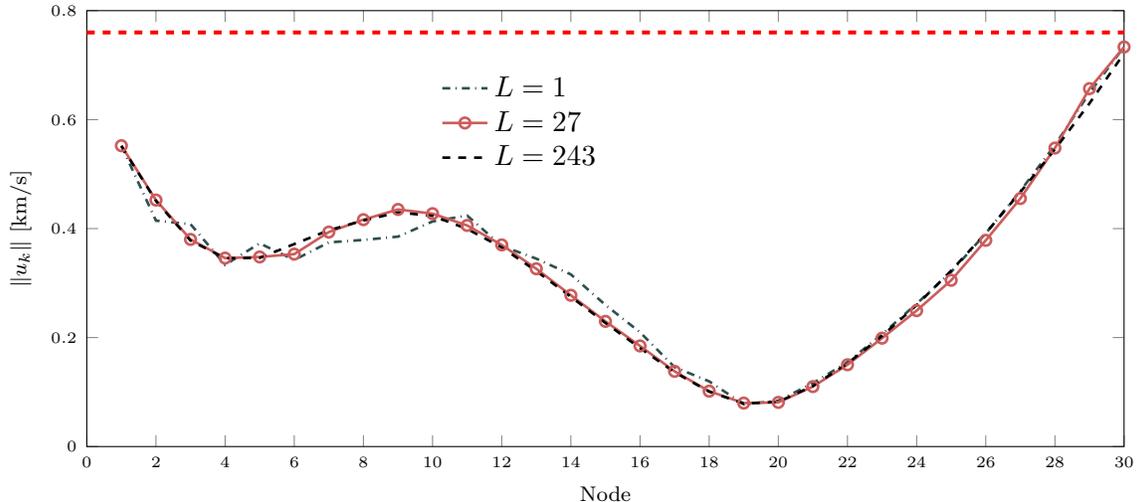
\begin{figure}[h!]
    \centering
    \setlength\figurewidth{14.5cm}
    \setlength\figureheight{5.8cm}
\definecolor{mycolor1}{rgb}{0.18359,0.30859,0.30859}%
\definecolor{mycolor2}{rgb}{0.80078,0.35938,0.35938}%
\begin{tikzpicture}

\begin{axis}[%
width=0.951\figurewidth,
height=\figureheight,
at={(0\figurewidth,0\figureheight)},
scale only axis,
xmin=0,
xmax=30,
xlabel style={font=\color{white!15!black}},
xlabel={Node},
ymin=0,
ymax=0.8,
ylabel style={font=\color{white!15!black}},
ylabel={$\| u_k \|$ [km/s]},
axis background/.style={fill=white},
legend style={legend cell align=left, align=left, draw=none, fill=none, anchor=north west, at={(axis cs:10, 0.7)}},
ylabel style={font=\color{white!15!black}, font=\scriptsize},
xlabel style={font=\color{white!15!black}, font=\scriptsize}
]
\addplot [color=mycolor1, dashdotted, line width=1.0pt]
  table[row sep=crcr]{%
1	0.552005187488678\\
2	0.414599386353313\\
3	0.408094032008194\\
4	0.333597336734968\\
5	0.372941372846277\\
6	0.342493499194859\\
7	0.374680080617406\\
8	0.37923380581514\\
9	0.385239248068387\\
10	0.413118881503102\\
11	0.423894660235579\\
12	0.368242038643512\\
13	0.34488406243814\\
14	0.315949676097418\\
15	0.259889354306328\\
16	0.209805377067699\\
17	0.145746335604464\\
18	0.119633328514823\\
19	0.0774772739165712\\
20	0.084162498790576\\
21	0.11631448966887\\
22	0.154806328761975\\
23	0.206011583265266\\
24	0.262069949885561\\
25	0.320473365671205\\
26	0.391249267972924\\
27	0.469005488969096\\
28	0.553116266075651\\
29	0.649188041164514\\
30	0.733552015428021\\
};
\addlegendentry{$L = 1$}

\addplot [color=mycolor2, line width=1.0pt, mark=o, mark options={solid, mycolor2}]
  table[row sep=crcr]{%
1	0.552147811236718\\
2	0.452400111041786\\
3	0.380117301843089\\
4	0.345908264167036\\
5	0.347963150549404\\
6	0.353177363097167\\
7	0.393629401890042\\
8	0.416398544198607\\
9	0.435139987167702\\
10	0.42741728058807\\
11	0.405805785400027\\
12	0.370099662292787\\
13	0.326342475385374\\
14	0.277609702282703\\
15	0.229737669748054\\
16	0.184683897394294\\
17	0.138087378339755\\
18	0.101653905676553\\
19	0.0797601118311775\\
20	0.0810408549070836\\
21	0.110027976985801\\
22	0.150131508959213\\
23	0.198902504466833\\
24	0.249469742082034\\
25	0.305176545469685\\
26	0.378474856790362\\
27	0.455185502206977\\
28	0.547793882856681\\
29	0.656897840157111\\
30	0.733320002910084\\
};
\addlegendentry{$L = 27$}

\addplot [color=black, dashed, line width=1.0pt]
  table[row sep=crcr]{%
1	0.551739906790424\\
2	0.450771444763824\\
3	0.378082020645298\\
4	0.345651898153761\\
5	0.346266588544189\\
7	0.396991578031038\\
8	0.414827771443715\\
9	0.429667101313434\\
10	0.423463820893808\\
11	0.399381856581329\\
12	0.365499839884215\\
13	0.322211648714465\\
14	0.274863540845459\\
15	0.226859656388775\\
16	0.180770004270979\\
17	0.136367636866307\\
18	0.10079222852108\\
19	0.0789822277103198\\
20	0.0821553327881759\\
21	0.111099086576132\\
22	0.152541775116443\\
23	0.203005273059457\\
24	0.259597540680975\\
25	0.322819000901085\\
26	0.39168621540599\\
27	0.467003148242174\\
28	0.545428936735707\\
29	0.629860504397392\\
30	0.721695630018321\\
};
\addlegendentry{$L = 243$}

\addplot [color=red, dashed, line width=1.5pt, forget plot]
  table[row sep=crcr]{%
0	0.76\\
30	0.76\\
};
\end{axis}
\end{tikzpicture}%
    \caption{Comparison of average closed-loop control magnitude for each scenario. Maximum allowable $\bm{\Delta V}$ is indicated by the red dashed line.}
    \label{fig::uCLcompare}
\end{figure}

Lastly, a comparison of the average (across MC trials) closed-loop control magnitude for each scenario is presented. The average magnitude history is illustrated in \bfig{fig::uCLcompare}. From inspection, the behavior of the average is similar for all scenarios, with the single Gaussian case being slightly more distinguished from the others. A measure of control cost in the Monte Carlo sense is computed by summing the average magnitude across the nodes. The result is an average total $\Delta V$ usage for each scenario and is given in \tbl{tbl::sumDV}. The GM model scenarios are able to beat the single Gaussian scenario by 100 m/s. This is directly a product of the GM-based ``solved for" mean states (used in \bequ{eqn::control}) being closer to the true state pdf mean, which results in a lower feedback control term than the single Gaussian model. This highlights the usefulness of proper uncertainty characterization in not only accuracy, but also in cost. 
\begin{table}[h!]
\centering
\caption{Sum of average closed-loop control magnitudes across nodes for each scenario.}
\label{tbl::sumDV}
    \begin{NiceTabular}{cc} \hline \hline
        $L$ & Total $\Delta V$ [km/s] \\
        \hline 
        1 & 10.181 \\
        27 & 10.084  \\
        243 & 10.081 \\
        \hline \hline
    \end{NiceTabular}
\end{table}

\section*{Conclusions}
In this paper, a Gaussian mixture formulation for solving a chance-constrained, stochastic optimal control problem, such as frequently encountered in finite-thrust spacecraft trajectory optimization, is presented. The problem statement is to find a control law that steers a stochastic, nonlinear system from one distribution to another whilst minimizing an objective and adhering to control chance constraints, otherwise known as distribution steering. The problem is decomposed into a recursive hybrid approach that uses a Gaussian mixture approximation of the state distribution for propagation, and a Gaussian approximation for the control execution. The result is a nonlinear, nonconvex program that can be used to generate feasible and locally optimal solutions. An Earth-to-Mars impulsive transfer problem is used to analyze the performance of the optimization. Additionally, a Monte Carlo analysis is conducted to validate the Gaussian mixture methodology and showcase its advantages in comparison to using a single Gaussian model.

For the Earth-to-Mars transfer, the optimization problem was solved for a single Gaussian case and two Gaussian mixture cases. Although constraints are met for all cases during the optimization step, Monte Carlo analysis with 50,000 samples shows that the single Gaussian case suffered the most from violation of the terminal covariance constraint. Conversely, the Gaussian mixture cases significantly lower the degree of closed-loop constraint violation. The significance of this result is that even though the optimization constraints were formulated assuming a single Gaussian, which breaks down in the face of nonlinear dynamics, embedding the Gaussian mixture methodology into the propagation alleviates these modeling errors. The Gaussian mixture cases also have lower closed-loop control cost compared to the single Gaussian case. All cases satisfy the control magnitude chance constraint as shown in Monte Carlo analysis. However, the current formulation of the chance constrained uncertainty steering problem does not model the control law or the chance constraints with a rigorous Gaussian mixture embedding. Constructing a fully mature Gaussian mixture version without at-node Gaussian approximations, as well as addressing additional, more complex, trajectory optimization problems are avenues of future effort.

\bibliographystyle{ieeetr}   % Number the references.
\bibliography{main}   % Use references.bib to resolve the labels.

\begin{thebibliography}{10}

\bibitem{bryson2018applied}
A.~E. Bryson and Y.-C. Ho, {\em Applied optimal control: optimization,
  estimation and control}.
\newblock Wiley, 1975.

\bibitem{hull1997conversion}
D.~G. Hull, ``Conversion of optimal control problems into parameter
  optimization problems,'' {\em Journal of guidance, control, and dynamics},
  vol.~20, no.~1, pp.~57--60, 1997.

\bibitem{hotz1987covariance}
A.~Hotz and R.~E. Skelton, ``Covariance control theory,'' {\em International
  Journal of Control}, vol.~46, no.~1, pp.~13--32, 1987.

\bibitem{bakolas2016optimal}
E.~Bakolas, ``Optimal covariance control for discrete-time stochastic linear
  systems subject to constraints,'' in {\em 2016 IEEE 55th Conference on
  Decision and Control (CDC)}, pp.~1153--1158, IEEE, 2016.

\bibitem{wainwright2019high}
M.~J. Wainwright, {\em High-dimensional statistics: A non-asymptotic
  viewpoint}, vol.~48.
\newblock Cambridge university press, 2019.

\bibitem{nemirovski2012safe}
A.~Nemirovski, ``On safe tractable approximations of chance constraints,'' {\em
  European Journal of Operational Research}, vol.~219, no.~3, pp.~707--718,
  2012.

\bibitem{pagnoncelli2009sample}
B.~K. Pagnoncelli, S.~Ahmed, and A.~Shapiro, ``Sample average approximation
  method for chance constrained programming: theory and applications,'' {\em
  Journal of optimization theory and applications}, vol.~142, no.~2,
  pp.~399--416, 2009.

\bibitem{geletu2015tractable}
A.~Geletu, M.~Kl{\"o}ppel, A.~Hoffmann, and P.~Li, ``A tractable approximation
  of non-convex chance constrained optimization with non-gaussian
  uncertainties,'' {\em Engineering Optimization}, vol.~47, no.~4,
  pp.~495--520, 2015.

\bibitem{balci2023density}
I.~Balci and E.~Bakolas, ``Density steering of gaussian mixture models for
  discrete-time linear systems,'' {\em arXiv preprint arXiv:2311.08500}, 2023.

\bibitem{boone2022spacecraft}
S.~Boone and J.~McMahon, ``Spacecraft maneuver design with non-gaussian chance
  constraints using gaussian mixtures,'' in {\em 2022 AAS/AIAA Astrodynamics
  Specialist Conference}, 2022.

\bibitem{marmo2022chance}
N.~Marmo, A.~Zavoli, {\em et~al.}, ``Chance-constraint optimization of
  interplanetary trajectories with a hybrid multiple-shooting approach,'' in
  {\em Proceedings of the 2022 AAS/AIAA Astrodynamics Specialist Conference},
  2022.

\bibitem{ross2004find}
I.~Ross, ``How to find minimum-fuel controllers,'' in {\em AIAA Guidance,
  Navigation, and Control Conference and Exhibit}, p.~5346, 2004.

\bibitem{skaf2010design}
J.~Skaf and S.~P. Boyd, ``Design of affine controllers via convex
  optimization,'' {\em IEEE Transactions on Automatic Control}, vol.~55,
  no.~11, pp.~2476--2487, 2010.

\bibitem{bertsekas1978stochastic}
D.~P. Bertsekas and S.~E. Shreve, {\em Stochastic Optimal Control: The
  Discrete-Time Case}.
\newblock USA: Academic Press, Inc., 1978.

\bibitem{sarkka2023bayesian}
S.~S{\"a}rkk{\"a} and L.~Svensson, {\em Bayesian filtering and smoothing},
  vol.~17.
\newblock Cambridge university press, 2023.

\bibitem{mclachlan2019finite}
G.~J. McLachlan, S.~X. Lee, and S.~I. Rathnayake, ``Finite mixture models,''
  {\em Annual review of statistics and its application}, vol.~6, pp.~355--378,
  2019.

\bibitem{li2009novel}
Y.~Li and L.~Li, ``A novel split and merge em algorithm for gaussian mixture
  model,'' in {\em 2009 Fifth International Conference on Natural Computation},
  vol.~6, pp.~479--483, IEEE, 2009.

\bibitem{hanebeck2003progressive}
U.~D. Hanebeck, K.~Briechle, and A.~Rauh, ``Progressive bayes: a new framework
  for nonlinear state estimation,'' in {\em Multisensor, Multisource
  Information Fusion: Architectures, Algorithms, and Applications 2003},
  vol.~5099, pp.~256--267, SPIE, 2003.

\bibitem{demars2013entropy}
K.~J. DeMars, R.~H. Bishop, and M.~K. Jah, ``Entropy-based approach for
  uncertainty propagation of nonlinear dynamical systems,'' {\em Journal of
  Guidance, Control, and Dynamics}, vol.~36, no.~4, pp.~1047--1057, 2013.

\bibitem{sarkka2019applied}
S.~S{\"a}rkk{\"a} and A.~Solin, {\em Applied stochastic differential
  equations}, vol.~10.
\newblock Cambridge University Press, 2019.

\bibitem{gates1963simplified}
C.~R. Gates, ``A simplified model of midcourse maneuver execution errors,''
  tech. rep., Jet Propulsion Lab, 1963.

\bibitem{sun2023hybrid}
P.~Sun, C.~Colombo, M.~Trisolini, and S.~Li, ``Hybrid gaussian mixture
  splitting techniques for uncertainty propagation in nonlinear dynamics,''
  {\em Journal of Guidance, Control, and Dynamics}, vol.~46, no.~4,
  pp.~770--780, 2023.

\bibitem{ridderhof2019nonlinear}
J.~Ridderhof, K.~Okamoto, and P.~Tsiotras, ``Nonlinear uncertainty control with
  iterative covariance steering,'' in {\em 2019 IEEE 58th Conference on
  Decision and Control (CDC)}, pp.~3484--3490, IEEE, 2019.

\bibitem{ridderhof2020chance}
J.~Ridderhof, J.~Pilipovsky, and P.~Tsiotras, ``Chance-constrained covariance
  control for low-thrust minimum-fuel trajectory optimization,'' in {\em 2020
  AAS/AIAA Astrodynamics Specialist Conference}, pp.~9--13, 2020.

\bibitem{boggs1995sequential}
P.~T. Boggs and J.~W. Tolle, ``Sequential quadratic programming,'' {\em Acta
  numerica}, vol.~4, pp.~1--51, 1995.

\bibitem{ozaki2020tube}
N.~Ozaki, S.~Campagnola, and R.~Funase, ``Tube stochastic optimal control for
  nonlinear constrained trajectory optimization problems,'' {\em Journal of
  Guidance, Control, and Dynamics}, vol.~43, no.~4, pp.~645--655, 2020.

\bibitem{gill2005snopt}
P.~E. Gill, W.~Murray, and M.~A. Saunders, ``Snopt: An sqp algorithm for
  large-scale constrained optimization,'' {\em SIAM review}, vol.~47, no.~1,
  pp.~99--131, 2005.

\end{thebibliography}

\end{document}